\documentclass[12pt]{article} \usepackage{amssymb}

\newtheorem{thm}{Theorem}[section]
\newtheorem{defn}[thm]{Definition}
\newtheorem{prop}[thm]{Proposition}

\newtheorem{rema}[thm]{Remark}

\newcommand{\pfbox}{\hspace*{\fill}\mbox{$\square$}}
\newcommand{\Z}{\mathbb{Z}_+}

\begin{document}

\begin{center}
\begin{Large}
{\bf  $N=1$ Neveu-Schwarz vertex operator superalgebras over 
Grassmann algebras and with odd formal variables}
\end{Large}

\medskip

Katrina Barron\footnote{Supported in part by an AAUW American
Dissertation Fellowship and by a University of California President's
Postdoctoral Fellowship}\\ Department of Mathematics, University of
California, Santa Cruz, CA 95064\\
\end{center}

\hspace{1.5cm}
\begin{abstract}
The notions of $N = 1$ Neveu-Schwarz vertex operator superalgebra 
over a Grassmann algebra and with odd formal variables and of 
$N = 1$ Neveu-Schwarz vertex operator superalgebra over a  
Grassmann algebra and without odd formal variables are introduced, 
and we show that the respective categories of such objects are 
isomorphic.  The weak supercommutativity and weak associativity 
properties for an $N = 1$ Neveu-Schwarz vertex operator 
superalgebra with odd formal variables are established, and we 
show that in the presence of the other axioms, weak 
supercommutativity and weak associativity are equivalent to the 
Jacobi identity.  In addition, we prove the supercommutativity 
and associativity properties for an $N = 1$  Neveu-Schwarz vertex 
operator superalgebra with odd formal variables and show that in 
the presence of the other axioms, supercommutativity and 
associativity are equivalent to the Jacobi identity. 
\end{abstract}

\section{Introduction}

In this paper, we introduce the notion of {\it $N = 1$ Neveu-Schwarz 
vertex operator superalgebra over a Grassmann algebra and with odd 
formal variables} and study the consequences of this notion.  We also 
introduce the notion of {\it $N = 1$ Neveu-Schwarz vertex operator 
superalgebra over a Grassmann algebra and without odd formal 
variables} and show that these two notions are equivalent in that the
corresponding categories of such objects are isomorphic.  Though many 
more or less equivalent notions of vertex operator superalgebra have  
been formulated (cf. \cite{T}, \cite{G}, \cite{FFR}, \cite{DL}, and 
\cite{KW}), extending the notion of vertex operator algebra \cite{Bo},
\cite{FLM}, \cite{FHL}, none of these notions give the most natural 
setting for the corresponding supergeometry of $N=1$ superconformal 
field theory (cf. \cite{Fd}, \cite{D}).  Some lack a representation 
of the full $N=1$ Neveu-Schwarz algebra rather than just the 
Virasoro algebra and none are defined over a general Grassmann 
algebra.    However, if one extends, for example, the notion of 
``$N=1$ NS-type SVOA" as given in \cite{KW} to be a module over a 
general Grassmann algebra instead of just a vector space over 
$\mathbb{C}$, then one does obtain an equivalent notion to ours.  
Since the notions of $N=1$ Neveu-Schwarz vertex operator 
superalgebra over a Grassmann algebra and with and without odd 
formal variables are equivalent, one might very well ask why we 
should want to complicate matters by adding the odd formal 
variables.   Below we give some motivation for including these odd 
components as  well as motivation for working over a Grassmann 
algebra and requiring  a representation of the $N=1$ Neveu-Schwarz 
algebra.  

In \cite{B thesis}, we give a rigorous foundation to the
correspondence between the geometric and algebraic aspects of 
genus-zero holomorphic $N=1$ superconformal field theory (cf. 
\cite{Fd}), following the work of Huang \cite{H book}, \cite{H thesis} 
by, introducing the notion of $N=1$ supergeometric vertex operator 
superalgebra and proving that the category of such objects is 
isomorphic to the category of $N = 1$ Neveu-Schwarz vertex operator 
superalgebras over a Grassmann algebra and with (or without) odd formal 
variables.  The supergeometry of $N=1$ superconformal field theory
is defined over a Grassmann algebra (cf. \cite{D}, \cite{R}), which is  
why it is more natural to work over such an algebra rather than just 
$\mathbb{C}$.  This supergeometry and the notion of $N=1$ supergeometric 
vertex operator superalgebra involve the moduli space of $N=1$ 
superspheres with tubes (corresponding to the interaction of incoming
and outgoing superstrings) which naturally has even and odd variables
in the underlying Grassmann algebra. 

Recall that in a vertex operator algebra, the Virasoro element $L(-1)$
is related to the differential operator $\frac{\partial}{\partial
z}$ via the $L(-1)$-derivative property (see \cite{FHL}).  This 
correspondence can be thought of as a correspondence between the 
algebraic setting and the geometric setting of genus-zero holomorphic
conformal field theory as developed in \cite{H thesis} and \cite{H 
book} following \cite{BPZ}, \cite{FS}, \cite{V}, \cite{S}.  In the 
geometry of conformal field theory, local coordinates are conformal in 
that they transform the differential operator $\frac{\partial}{\partial 
z}$ homogeneously of degree one for $z$ a complex variable.  In $N=1$ 
superconformal field theory, the natural setting is $N=1$ 
supergeometry, and local coordinates which transform a certain
superdifferential operator $D$ homogeneously of degree one where
$D$ satisfies $D^2 = \frac{\partial}{\partial z}$.  Such an operator 
is given by $D = \frac{\partial}{\partial \theta} + 
\theta \frac{\partial}{\partial z}$ where $z$ is an even variable 
over a Grassmann algebra and $\theta$ is an odd variable over a 
Grassmann algebra (cf. \cite{Fd}, \cite{B thesis}, \cite{B thesis 
announcement}).  In addition, in \cite{B thesis},  we show that the 
Lie superalgebra of infinitesimal local coordinate transformations 
for an $N=1$ supersphere with tubes is isomorphic to the $N=1$ 
Neveu-Schwarz algebra \cite{NS} with central charge zero.

Thus in considering what should be the corresponding
superalgebraic setting for $N=1$ supergeometric vertex operator
superalgebra, we naturally want to have a representation of the $N=1$ 
Neveu-Schwarz algebra in which the element $G(-\frac{1}{2})$ has the 
supercommutator $\frac{1}{2}[G(-\frac{1}{2}), G(-\frac{1}{2})] = L(-1)$ 
(which corresponds to $G(-\frac{1}{2})^2 = L(-1)$ in the universal 
enveloping algebra of the Neveu-Schwarz algebra).  Furthermore, this
operator $G(-\frac{1}{2})$ should be related to the 
superdifferential operator $D = \frac{\partial}{\partial \theta} + 
\theta \frac{\partial}{\partial z}$ consistent with the relationship 
between the $L(-1)$ operator and $\frac{\partial}{\partial z}$.  In 
our notion of $N=1$ Neveu-Schwarz vertex operator superalgebra over a 
Grassmann algebra and with odd formal variables, we have an $N=1$ 
Neveu-Schwarz algebra element giving rise to a  representation of the 
$N=1$ Neveu-Schwarz algebra with a corresponding 
$G(-\frac{1}{2})$-derivative property (property 
(\ref{G(-1/2)-derivative}) in Definition \ref{VOSA definition} below) 
giving the explicit connection between the endomorphism $G(-\frac{1}{2})$ 
and the superderivation $D$.  The $L(-1)$-derivative property then 
follows {}from this more fundamental property.

Furthermore, the correlation functions of $N=1$ superconformal field 
theory (cf. \cite{Fd}) are superanalytic superfunctions in one even and 
one odd variable and, in the genus-zero holomorphic case, correspond to 
the correlation functions arising {}from an $N=1$ Neveu-Schwarz vertex 
operator superalgebra over a Grassmann algebra and with odd formal 
variables as studied in Section 9 of this paper, i.e., our inclusion of 
the odd formal variables gives the entire correlation function 
explicitly.

Thus the fact that we are working over a general Grassmann algebra and
the presence of a representation of the $N=1$ Neveu-Schwarz algebra, 
the $G(-\frac{1}{2})$-derivative property, and the odd formal variable 
component of the vertex operators in our definition of $N=1$  
Neveu-Schwarz vertex operator algebra over a Grassmann algebra and with  
odd formal variables give a more natural correspondence to the  
supergeometry of $N=1$ superconformal field theory.

For the remainder of this paper we will often use the term ``vertex 
operator superalgebra" to mean ``$N=1$  Neveu-Schwarz vertex operator 
superalgebra over a Grassmann algebra" unless otherwise noted.

That the notions of vertex operator superalgebra with and without odd 
formal variables are equivalent is due to the fact that all of the 
information for the odd variable components of the vertex operators
is contained in the representation of the Neveu-Schwarz algebra.  In 
fact, given a vertex operator superalgebra without odd formal variables
there is only one way to form the odd variable component in order
to obtain a vertex operator superalgebra with odd formal variables and
that is by employing the $G(-\frac{1}{2})$ representative element (see
consequence (\ref{odd relationship}) of Definition \ref{VOSA definition} 
below).

After introducing the notions of vertex operator superalgebra with and
without odd formal variables and proving their equivalence, we 
formulate the properties of weak supercommutativity and weak 
associativity for a vertex operator superalgebra with odd formal 
variables.  Following \cite{DL}, \cite{L1} and \cite{L2}, we show that 
in the presence of the other axioms for a vertex operator superalgebra 
with odd formal variables, weak supercommutativity and weak 
associativity are equivalent to the Jacobi identity.  Then following 
\cite{FLM} and \cite{FHL}, we look at expansions of rational 
superfunctions and formulate the properties of rationality of 
products and iterates, supercommutativity and associativity for a 
vertex operator superalgebra with odd formal variables.  Finally, 
following \cite{FHL}, we show that in the presence of the other 
axioms, rationality of products and iterates, supercommutativity and 
associativity are equivalent to the Jacobi identity.  These properties 
are used in the proof that the category of $N=1$ supergeometric vertex 
operator superalgebras is isomorphic to the category of vertex operator 
superalgebras with odd formal variables given in \cite{B thesis} and  
announced in \cite{B thesis announcement}.

The results in this paper follow the theory of vertex operator 
algebras as developed in \cite{Bo}, \cite{FLM}, \cite{FHL},
\cite{DL} and \cite{L2} and are contained in the author's Ph.D. 
dissertation \cite{B thesis}.

We would like to thank James Lepowsky and Yi-Zhi Huang for their
advice, support and expert comments during the writing of the 
dissertation {}from which this paper is derived, and for 
additional comments regarding the exposition of this paper.

\section{Grassmann algebras and the $N=1$ Neveu-Schwarz algebra}

Let $\mathbb{Z}_2$ denote the integers modulo two.  For a 
$\mathbb{Z}_2$-graded vector space $V = V^0 \oplus V^1$, define the 
{\it sign function} $\eta$ on the homogeneous subspaces of $V$ by 
$\eta(v) = i$ for $v \in V^i$, $i = 0,1$.  If $\eta(v) = 0$, we say 
that $v$ is {\it even}, and if $\eta(v) = 1$, we say that $v$ is 
{\it odd}. (Note that with this definition, the sign function is 
double-valued at zero.  However, in practice this is never an issue,
as can be seen in the definitions below.)

A {\it superalgebra} is an (associative) algebra $A$ (with identity $1
\in A$), such that
\begin{eqnarray*}
&(i)& \mbox{$A$ is a $\mathbb{Z}_2$-graded algebra}\\ &(ii)& \mbox{$ab
= (-1)^{\eta(a) \eta(b)} ba$ for $a,b$ homogeneous in $A$.}
\hspace{1.5in}
\end{eqnarray*}

The exterior algebra over a vector space $U$, denoted $\bigwedge(U)$, 
has the structure of a superalgebra.  Fix $U_L$ to be an 
$L$-dimensional vector space over $\mathbb{C}$ for $L \in \mathbb{N}$ 
such that $U_L \subset U_{L+1}$.  We denote $\bigwedge(U_L)$ by 
$\bigwedge_L$ and call this the {\it Grassmann algebra on $L$  
generators}.  Note that $\bigwedge_L \subset \bigwedge_{L+1}$, and 
taking the direct limit as $L \rightarrow \infty$, we have the
{\it infinite Grassmann algebra} denoted by $\bigwedge_\infty$.  We 
use the notation $\bigwedge_*$ to denote a Grassmann algebra, finite 
or infinite.  Since this paper is motivated by the geometry of $N=1$ 
superconformal field theory in which one considers superanalytic 
structures, for the purposes of this paper we have defined 
$\bigwedge_*$ over $\mathbb{C}$.  However, we could just as well 
have formulated the results that follow for Grassmann algebras over 
any field of characteristic zero.

A $\mathbb{Z}_2$-graded vector space $\mathfrak{g}$ is said to be a
{\it Lie superalgebra} if it has a bilinear operation $[\cdot,\cdot]$
such that for $u,v$ homogeneous in $\mathfrak{g}$,
\begin{eqnarray*}
&(i)& [u,v] \in \mathfrak{g}^{(\eta(u) + \eta(v)) \mathrm{mod} \; 2}\\ 
&(ii)& [u,v] = -(-1)^{\eta(u)\eta(v)}[v,u] \hspace{1.75in} 
\mbox{(skew-symmetry)}\\
&(iii)& (-1)^{\eta(u)\eta(w)}[[u,v],w] + (-1)^{\eta(v)\eta(u)}[[v,w],u] \\
& & \hspace{1.35in} + \; (-1)^{\eta(w)\eta(v)}[[w,u],v] = 0. \qquad 
\mbox{(Jacobi identity)}
\end{eqnarray*}

The Virasoro algebra is the Lie algebra with central charge $d$,
basis consisting of the central element $d$ and $L_n$, for $n \in 
\mathbb{Z}$, and commutation relations  
\begin{equation}\label{Virasoro relation}
[L_m ,L_n] = (m - n)L_{m + n} + \frac{1}{12} (m^3 - m) \delta_{m + n 
, 0} \; d ,
\end{equation}
for $m, n \in \mathbb{Z}$.  The $N=1$ Neveu-Schwarz Lie superalgebra
is a super-extension of the Virasoro algebra by the odd elements 
$G_{n + \frac{1}{2}}$, for $n \in \mathbb{Z}$.  That is, the $N=1$
Neveu-Schwarz algebra has a basis consisting of the central element 
$d$, $L_n$ and $G_{n + \frac{1}{2}}$, for $n \in \mathbb{Z}$, with 
supercommutation relations 
\begin{eqnarray}
\left[ G_{m + \frac{1}{2}},L_n \right] &=& \Bigl(m - \frac{n - 1}{2} \Bigr) G_{m
+ n + \frac{1}{2}}  \label{Neveu-Schwarz relation1} \\  
\left[ G_{m + \frac{1}{2}} , G_{n - \frac{1}{2}} \right] &=& 2L_{m +
n} + \frac{1}{3} (m^2 + m) \delta_{m + n , 0} \; d \label{Neveu-Schwarz relation2}
\end{eqnarray}
in addition to (\ref{Virasoro relation}).

\section{Delta functions with odd formal variables}

Let $x_1$ and $x_2$ be formal variables which commute with each other
and with $\bigwedge_*$.  We call such a variable an {\it even formal
variable}.  Let $\varphi_1$ and $\varphi_2$ be formal variables which
commute with $x_1$, $x_2$ and $\bigwedge_*^0$ and anticommute with each
other and $\bigwedge_*^1$.  We call such variables {\it odd formal
variables}.  Note that the square of any odd formal variable is zero.  
Thus for any formal Laurent series $f(x) \in  \bigwedge_*[[x,x^{-1}]]$, 
we can define 
\begin{equation}\label{expansion}
f(x + \varphi_1 \varphi_2) = f(x) + \varphi_1
\varphi_2 f'(x) \; \in \mbox{$\bigwedge_*$} [[x,x^{-1}]][\varphi_1][\varphi_2]. 
\end{equation}

In order to formulate the notion of vertex operator superalgebra with
odd formal variables, we would like to use the {\it formal 
$\delta$-function at $x=1$} given by 
\[\delta(x) = \sum_{n \in \mathbb{Z}} x^n ,\] 
which is a fundamental ingredient in the formal calculus underlying the 
theory of vertex operator algebras.  However, we will want to extend 
the $\delta$-function to be defined for certain expressions involving
both even and odd formal variables using (\ref{expansion}). 

Let $x$, $x_0$, $x_1$, and $x_2$ be even formal variables.  Following
the treatment of formal calculus for even formal variables given in 
\cite{FHL}, we have  
\begin{equation}\label{delta 1}
f(x) \delta (x) = f(1) \delta (x) \quad \mbox{for} \; \; f(x) \in
\mathbb{C}[x,x^{-1}] .
\end{equation}
Let $V$ be a vector space over $\mathbb{C}$.  For $X(x_1,x_2) \in 
(\mbox{End} \; V)[[x_1, x_1^{-1},x_2, x_2^{-1}]]$ such that
\[\lim_{x_1 \rightarrow x_2} X(x_1,x_2) \quad \mbox{exists, i.e.,}
\quad X(x_1,x_2)|_{x_1 = x_2} \quad \mbox{exists} \]
(that is, when $X(x_1,x_2)$ is applied to any element of $V$, setting
the variables equal leads to only finite sums in $V$), we have
\begin{equation}\label{delta multiplication}
X(x_1,x_2) \delta \biggl( \frac{x_1}{x_2} \biggr) = X(x_2,x_2) \delta
\biggl( \frac{x_1}{x_2} \biggr) .
\end{equation}
Let $\mathbb{N}$ denote the natural numbers.  For the three-variable 
generating function
\[\delta \biggl( \frac{x_1 - x_2}{x_0} \biggr) = \sum_{n \in \mathbb{Z}}
\frac{(x_1 - x_2)^n}{x_0^n} = \sum_{n \in \mathbb{Z}, m \in \mathbb{N}}
(-1)^m \biggl( \! \! \begin{footnotesize} 
\begin{array}{c}   
n \\ 
m  \\
\end{array} \end{footnotesize} \! \! \biggr)
x_0^{-n} x_1^{n-m} x_2^m ,\]  
we have
\begin{equation}\label{delta 2 terms}
x_1^{-1} \delta \biggl( \frac{x_2 + x_0}{x_1} \biggr) = x_2^{-1}
\delta \biggl( \frac{x_1 - x_0}{x_2} \biggr) 
\end{equation}
and
\begin{equation}\label{delta 3 terms}
x_0^{-1} \delta \biggl( \frac{x_1 - x_2}{x_0} \biggr) - x_0^{-1}
\delta \biggl( \frac{x_2 - x_1}{-x_0} \biggr) = x_2^{-1} \delta \biggl(
\frac{x_1 - x_0}{x_2} \biggr) ,
\end{equation}
where we use the convention that any expression such as 
$(x_1 - x_2)^n$ for $n \in \mathbb{Z}$, is understood to be expanded 
in positive powers of $x_2$.  We will continue to use this convention
throughout the rest of this work.  

Notice that in the spirit of the $\delta$-function multiplication
principle (\ref{delta multiplication}), the expressions on both sides
of (\ref{delta 2 terms}) and the three terms occurring in (\ref{delta
3 terms}) all correspond to the same formal substitution 
\[x_0 = x_1 - x_2 . \] 

Let $\varphi$, $\varphi_1$, and $\varphi_2$ be odd formal variables.  
Extending the above results, we see that for any $f(x,\varphi) \in
\bigwedge_*[x,x^{-1},\varphi]$, the following is an immediate 
consequences of (\ref{delta 1}):
\[f(x,\varphi)\delta(x) = f(1,\varphi)\delta(x)  .\]
Let $V$ be a $\bigwedge_*$-module. If
$X(x_1,\varphi_1,x_2,\varphi_2) \in (\mbox{End} \;
V)[[x_1,x_1^{-1},x_2, x_2^{-1}]][\varphi_1, \varphi_2]$ such that  
\[\lim_{x_1 \rightarrow x_2} X(x_1,\varphi_1,x_2,\varphi_2) \quad
\mbox{exists, i.e.,} \quad X(x_1,\varphi_1,x_2,\varphi_2)|_{x_1 = x_2} 
\quad \mbox{exists} ,\]  
then we have the following immediate consequence of (\ref{delta
multiplication}):
\begin{equation}\label{delta multiplication with phis}
X(x_1,\varphi_1,x_2,\varphi_2) \delta \biggl( \frac{x_1}{x_2} \biggr) =
X(x_2,\varphi_1,x_2,\varphi_2) \delta \biggl(\frac{x_1}{x_2} \biggr) . 
\end{equation}

We have the following $\delta$-function of expressions involving three
even variables and two odd variables
\begin{eqnarray*}
\delta \biggl( \frac{x_1 - x_2 - \varphi_1 \varphi_2}{x_0} \biggr) &=&
\sum_{n \in \mathbb{Z}} (x_1 - x_2 - \varphi_1 \varphi_2)^n x_0^{-n} \\
&=& \sum_{n \in \mathbb{Z}} \biggl( (x_1 - x_2)^n - n
\varphi_1 \varphi_2 (x_1 - x_2)^{n - 1} \biggr) x_0^{-n} \\
&=& \delta \biggl( \frac{x_1 - x_2}{x_0} \biggr)  -
\varphi_1 \varphi_2 x_0^{-1} \delta' \biggl( \frac{x_1 - x_2}{x_0}
\biggr)
\end{eqnarray*} 
where 
\[\delta'(x) = \frac{d}{dx} \delta (x) = \sum_{n \in \mathbb{Z}} n
x^{n-1} ,\] 
and we use the convention that a function of even and odd variables
should be expanded about the even variables.  Conceptually, however,
we can also think of $\delta \left( \frac{x_1 - x_2 - \varphi_1
\varphi_2}{x_0} \right)$ as being a function in four even variables
where $\varphi_1 \varphi_2$ is considered as another even formal
variable $x_3$ with the property that $x_3^2 = 0$.  

Taking $\frac{\partial}{\partial x_0}$ of both sides of (\ref{delta 2
terms}), we obtain the following identity:
\begin{equation}\label{delta derivative 2 terms}
x_1^{-2} \delta' \biggl( \frac{x_2 + x_0}{x_1} \biggr) = - x_2^{-2}
\delta' \biggl( \frac{x_1 - x_0}{x_2} \biggr) .
\end{equation}
Taking $\frac{\partial}{\partial x_1}$ of both sides of (\ref{delta 3
terms}), we obtain:
\begin{equation}\label{delta derivative 3 terms}
x_0^{-2} \delta' \biggl( \frac{x_1 - x_2}{x_0} \biggr) - x_0^{-2}
\delta' \biggl( \frac{x_2 - x_1}{-x_0} \biggr) = x_2^{-2} \delta'
\biggl( \frac{x_1 - x_0}{x_2} \biggr) .
\end{equation}
        
Thus {}from (\ref{delta 2 terms}), (\ref{delta 3 terms}),
(\ref{delta derivative 2 terms}) and (\ref{delta derivative 3 terms}), 
we have 
\begin{equation}\label{delta 2 terms with phis}
x_1^{-1} \delta \biggl( \frac{x_2 + x_0 + \varphi_1 \varphi_2}{x_1}
\biggr) = x_2^{-1} \delta \biggl( \frac{x_1 - x_0 - \varphi_1
\varphi_2}{x_2} \biggr)  
\end{equation}
and
\begin{eqnarray}
x_0^{-1} \delta \biggl( \frac{x_1 - x_2 - \varphi_1 \varphi_2}{x_0}
\biggr) - x_0^{-1} 
\delta \biggl( \frac{x_2 - x_1 + \varphi_1 \varphi_2}{-x_0} \biggr) = 
\label{delta 3 terms with phis} \hspace{1.2in} \\
\hspace{2in} x_2^{-1} \delta \biggl( \frac{x_1 - x_0 - \varphi_1 
\varphi_2}{x_2} \biggr) . \nonumber
\end{eqnarray}
Notice that in the spirit of the $\delta$-function multiplication
principle (\ref{delta multiplication with phis}), the expressions on
both sides of (\ref{delta 2 terms with phis}) and the three terms
occurring in (\ref{delta 3 terms with phis}) all correspond to the
same formal substitution 
\[x_0 = x_1 - x_2  - \varphi_1 \varphi_2 . \]  
And thus, not surprisingly, for
\[X(x_0,\varphi_0,x_1,\varphi_1,x_2,\varphi_2) \in (\mbox{End} \;V)
[[x_0,x_0^{-1},x_1,x_1^{-1},x_2,x_2^{-1}]] [\varphi_1, \varphi_2] ,\] 
a formal substitution corresponding to $x_0 = x_1 - x_2 - \varphi_1
\varphi_2$ can be made as long as the resulting expression is
well defined, e.g., if 
\[X(x_1,\varphi_1,x_2,\varphi_2) \in (\mbox{End} \;V)
[[x_1,x_1^{-1}]]((x_2)) [\varphi_1, \varphi_2] ,\] 
then  
\begin{equation}\label{delta substitute}
\delta \biggl( \frac{x_2 + x_0 + \varphi_1 \varphi_2}{x_1} \biggr)
X(x_1,\varphi_1,x_2,\varphi_2) = \hspace{2.2in} 
\end{equation}
\[\hspace{1.5in} \delta \biggl( \frac{x_2 + x_0 + \varphi_1
\varphi_2}{x_1} \biggr) X(x_2 + x_0 + \varphi_1 \varphi_2,
\varphi_1,x_2,\varphi_2) .\]

The substitution $x_0 = x_1 - x_2 - \varphi_1 \varphi_2$ can be 
thought of as the even part of a superconformal shift of 
$(x_1, \varphi_1)$ by $(x_2, \varphi_2)$.  Formally, a power 
series $f(x_1,\varphi_1) \in \bigwedge_*[[x_1]][\varphi_1]$ is 
superconformal in $x_1$ and $\varphi_1$ if and only if $f$ 
satisfies $D\tilde{x} = \tilde{\varphi} D \tilde{\varphi}$ for 
$f(x_1,\varphi_1) = (\tilde{x}, \tilde{\varphi})$ and $D = 
\frac{\partial}{\partial \varphi_1} + \varphi_1 
\frac{\partial}{\partial x_1}$ (see \cite{B thesis}). (For a 
superanalytic function $f$ in one even variable and one odd variable
this condition is equivalent to requiring that $f$ transform the
super-differential operator $D$ homogeneously of degree one.)
Thus $f(x_1,\varphi_1) = (x_1 - x_2 - \varphi_1 \varphi_2, 
\varphi_1 - \varphi_2)$ is formally superconformal in $x_1$ and 
$\varphi_1$ since
\begin{eqnarray*}
D\tilde{x} &=& \Bigl(\frac{\partial}{\partial \varphi_1} + \varphi_1 
\frac{\partial}{\partial x_1}\Bigr) (x_1 - x_2 - \varphi_1 \varphi_2)\\
&=& -\varphi_2 + \varphi_1 \\
&=& (\varphi_1 - \varphi_2) \Bigl(\frac{\partial}{\partial \varphi_1} + \varphi_1 
\frac{\partial}{\partial x_1}\Bigr) (\varphi_1 -\varphi_2)\\
&=& \tilde{\varphi}D\tilde{\varphi}.
\end{eqnarray*}

\section{The notion of vertex operator superalgebra over $\bigwedge_*$
and with odd formal variables}

\begin{defn}\label{VOSA definition}
{\em A} ($N = 1$ Neveu-Schwarz) vertex operator superalgebra over 
$\bigwedge_*$ and with odd variables {\em is a $\frac{1}{2} 
\mathbb{Z}$-graded (by weight) $\bigwedge_*$-module which is 
also $\mathbb{Z}_2$-graded (by sign)  
\begin{equation}\label{vosa1}
V = \coprod_{n \in \frac{1}{2} \mathbb{Z}} V_{(n)} = \coprod_{n
\in \frac{1}{2}\mathbb{Z}} V_{(n)}^0 \oplus \coprod_{n \in
\frac{1}{2}\mathbb{Z}} V_{(n)}^1 = V^0 \oplus V^1 
\end{equation}  
such that 
\begin{equation}\label{vosa2}
\dim V_{(n)} < \infty \quad \mbox{for} \quad n \in \frac{1}{2}
\mathbb{Z} , 
\end{equation}
\begin{equation}\label{positive energy}
V_{(n)} = 0 \quad \mbox{for $n$ sufficiently small} , 
\end{equation}
equipped with a linear map $V \otimes V \longrightarrow V[[x,x^{-1}]]
[\varphi]$, or equivalently,
\begin{eqnarray*} 
V &\longrightarrow&  (\mbox{End} \; V)[[x,x^{-1}]][\varphi] \\
v  &\mapsto&  Y(v,(x,\varphi)) = \sum_{n \in \mathbb{Z}} v_n x^{-n-1} +
\varphi \sum_{n \in \mathbb{Z}} v_{n - \frac{1}{2}} x^{-n-1}
\end{eqnarray*}
where $v_n \in (\mbox{End} \; V)^{\eta(v)}$ and $v_{n -
\frac{1}{2}} \in (\mbox{End} \; V)^{(\eta(v) + 1)
\mbox{\begin{footnotesize} mod \end{footnotesize}} 2}$ for $v$ of 
homogeneous sign in $V$, $x$ is an even formal variable, and $\varphi$
is an odd formal variable, and where $Y(v,(x,\varphi))$ denotes the}
vertex operator associated with $v$, {\em and equipped also with two
distinguished homogeneous vectors $\mbox{\bf 1} \in V_{(0)}^0$ (the
{\em vacuum}) and $\tau \in V_{(\frac{3}{2})}^1$ (the {\em 
Neveu-Schwarz element}).  The following conditions are assumed for 
$u,v \in V$:   
\begin{equation}\label{truncation}
u_n v = 0 \quad \mbox{for $n \in \frac{1}{2} \mathbb{Z}$ sufficiently
large;} 
\end{equation}
\begin{equation}\label{vacuum identity}
Y(\mbox{\bf 1}, (x, \varphi)) = 1 \quad \mbox{(1 on the right being
the identity operator);} 
\end{equation}
the} creation property {\em holds:
\[Y(v,(x,\varphi)) \mbox{\bf 1} \in V[[x]][\varphi] \qquad \mbox{and}
\qquad \lim_{(x,\varphi) \rightarrow 0} Y(v,(x,\varphi)) \mathbf{1} =
v ; \] 
the} Jacobi identity {\em holds:  
\[x_0^{-1} \delta \biggl( \frac{x_1 - x_2 - \varphi_1 \varphi_2}{x_0}
\biggr) Y(u,(x_1, \varphi_1))Y(v,(x_2, \varphi_2)) \hspace{1.6in} \]
\[- (-1)^{\eta(u)\eta(v)} x_0^{-1} \delta \biggl( \frac{x_2 - x_1 + 
\varphi_1 \varphi_2}{-x_0} \biggr)Y(v,(x_2, \varphi_2))Y(u,(x_1,
\varphi_1)) \]
\[\hspace{1.1in} = \; x_2^{-1} \delta \biggl( \frac{x_1 - x_0 - \varphi_1
\varphi_2}{x_2} \biggr) Y(Y(u,(x_0, \varphi_1 - \varphi_2))v,(x_2,
\varphi_2)) , \]
for $u,v$ of homogeneous sign in $V$; the Neveu-Schwarz algebra 
relations hold:
\begin{eqnarray*}
\left[L(m),L(n) \right] \! &=& \!(m - n)L(m + n) + \frac{1}{12} (m^3 - m)
\delta_{m + n , 0} (\mbox{rank} \; V) , \\ \label{V9}
\biggl[ G(m + \frac{1}{2}),L(n) \biggr] \! &=& \! (m - \frac{n - 1}{2} ) G(m
+ n + \frac{1}{2}) ,\\ \label{V10}
\biggl[ G (m + \frac{1}{2} ) , G(n - \frac{1}{2} ) \biggr] \! &=& \! 2L(m +
n) + \frac{1}{3} (m^2 + m) \delta_{m + n , 0} (\mbox{rank} \; V) , \label{V11}
\end{eqnarray*}
for $m,n \in \mathbb{Z}$, where 
\[G(n + \frac{1}{2}) = \tau_{n + 1}, \qquad \mbox{and} \qquad 2L(n) =
\tau_{n + \frac{1}{2}} \qquad \mbox{for} \; n \in \mathbb{Z} , \]
i.e., 
\begin{equation}\label{stress tensor}
Y(\tau,(x,\varphi)) = \sum_{n \in \mathbb{Z}} G (n + \frac{1}{2}) x^{-
n - \frac{1}{2} - \frac{3}{2}} \; + \; 2 \varphi \sum_{n \in \mathbb{Z}}
L(n) x^{- n - 2} ,
\end{equation}
and $\mbox{rank} \; V \in \mathbb{C}$; 
\begin{equation}\label{grading for vosa with}
L(0)v = nv \quad \mbox{for} \quad n \in \frac{1}{2} \mathbb{Z} \quad
\mbox{and} \quad v \in V_{(n)}; 
\end{equation}
\begin{equation}\label{G(-1/2)-derivative}
\biggl( \frac{\partial}{\partial \varphi} + \varphi
\frac{\partial}{\partial x} \biggr) Y(v,(x,\varphi)) =  Y(G(- 
\frac{1}{2})v,(x,\varphi)) . 
\end{equation} }
\end{defn}

The vertex operator superalgebra just defined is denoted by
\[(V,Y(\cdot,(x,\varphi)),\mbox{\bf 1},\tau) \]  
or for simplicity by $V$.

For such a vertex operator algebra $V$ over $\bigwedge_*$ for
$\bigwedge_* = \bigwedge_0 = \mathbb{C}$, the $\mathbb{Z}_2$-grading
is given by
\[ V^0 = \coprod_{n\in \mathbb{Z}} V_{(n)} \qquad \qquad V^1 = 
\coprod_{n\in \mathbb{Z} + \frac{1}{2}} V_{(n)} .\]
Extending $V$ to a vertex operator superalgebra over a general
$\bigwedge_*$ by $\bigwedge_* \otimes V$ changes the  
$\mathbb{Z}_2$-grading via $(\bigwedge_* \otimes V)^0 = \bigwedge_*^0 
\otimes V^0 + \bigwedge_*^1 \otimes V^1$ and  $(\bigwedge_* 
\otimes V)^1 = \bigwedge_*^1 \otimes V^0 + \bigwedge_*^0 \otimes V^1$.

Since $\bigwedge_L \subset \bigwedge_\infty$ for $L\in \mathbb{N}$,
any vertex operator superalgebra $V$ over $\bigwedge_L$ can be extended 
to a vertex operator superalgebra over $\bigwedge_\infty$ by, for
instance, defining the action of $\bigwedge_\infty \smallsetminus 
\bigwedge_L$ on $V$ to be trivial, or by taking the 
$\bigwedge_\infty$-module induced by the $\bigwedge_L$-module $V$.

Extending the proofs for the analogous properties in the non-super case 
{}from \cite{FHL}, we have the following consequences of the definition
of vertex operator superalgebra with odd formal variables.  
\begin{eqnarray}
L(n) \mbox{\bf 1} = 0 , &\mbox{and}& G(n + \frac{1}{2})
\mbox{\bf 1} = 0 , \; \mbox{for} \; n \geq -1 ;\\
G(-\frac{3}{2}) \mbox{\bf 1} &=& \tau ;\\
L(0) \tau &=& \frac{3}{2} \tau ;\\
Y(v,(x,\varphi)) \mbox{\bf 1} &=& e^{xL(-1) + \varphi
G(-\frac{1}{2})}v  \nonumber \\
&=& e^{xL(-1)}v + \varphi G(-\frac{1}{2})e^{xL(-1)} v ;
\end{eqnarray}
there exists $\omega = \frac{1}{2}G(-\frac{1}{2}) \tau \in V_{(2)}$ 
such that 
\begin{equation}
Y(\omega,(x,\varphi)) = \sum_{n \in \mathbb{Z}} L(n) x^{-n-2} -
\frac{\varphi}{2} \sum_{n \in \mathbb{Z}} (n + 1) G(n -
\frac{1}{2}) x^{-n-2} ;
\end{equation}
\begin{equation}
\mbox{wt} \; v_n = \mbox{wt} \; v - n - 1 , 
\end{equation}
for $n \in \frac{1}{2} \mathbb{Z}$ and for $v \in V$ of homogeneous 
weight.

Applying the $G(-\frac{1}{2})$-derivative property (\ref{G(-1/2)-derivative})
twice, we obtain the $L(-1)$-derivative property
\begin{equation}\label{L(-1)-derivative} 
\frac{\partial}{\partial x} Y(v, (x,\varphi)) = Y(L(-1)v,(x,\varphi)).
\end{equation}

Using the other properties, we see that the $N=1$ Neveu-Schwarz algebra
supercommutation relations are equivalent to:
\begin{equation}
Y(\tau, (x,\varphi))\tau = \frac{2}{3}(\mbox{rank} \; V) \mathbf{1} x^{-3}
+ 2 \omega x^{-1} + \varphi(3\tau x^{-2} + 2L(-1) \tau x^{-1}) + y
\end{equation}
where $y \in V[[x]][\varphi]$.
 
Taking $\mbox{Res}_{x_0}$ of the Jacobi identity and using 
the $\delta$-function identity (\ref{delta 3 terms with phis}), we obtain 
the following supercommutator formula where $\mbox{Res}_{x_0}$ of a power 
series in $x_0$ is the coefficient of $x_0^{-1}$.
\begin{equation}\label{bracket relation for a vosa}
[ Y(u, (x_1,\varphi_1)), Y(v,(x_2,\varphi_2))] = \hspace{2.6in} 
\end{equation}
\[\hspace{.5in} \mbox{Res}_{x_0} x_2^{-1} \delta \biggl( \frac{x_1 - x_0
- \varphi_1 \varphi_2}{x_2} \biggr) Y(Y(u,(x_0, \varphi_1 -
\varphi_2))v,(x_2, \varphi_2)) . \] 
{}From the $G(-\frac{1}{2})$-derivative property (\ref{G(-1/2)-derivative}) 
and the supercommutator formula (\ref{bracket relation for a vosa}), we have
\begin{equation}\label{odd relationship}
Y(v,(x,\varphi)) = \sum_{n \in \mathbb{Z}} v_n x^{-n-1} +
\varphi \sum_{n \in \mathbb{Z}} [G(- \frac{1}{2}), v_n ] x^{-n-1} ,
\end{equation}
i.e., $v_{n - \frac{1}{2}} =  [G(- \frac{1}{2}), v_n]$.

Taking $\mbox{Res}_{x_0} \mbox{Res}_{x_1}$ of the Jacobi identity, 
we find that for $u,v \in V$
\begin{equation}
\bigl[u_{-\frac{1}{2}}, Y(v,(x,\varphi)) \bigr] = Y(u_{-\frac{1}{2}}v,(x,\varphi))
\end{equation}
and in particular,
\begin{equation}
\bigl[u_{-\frac{1}{2}}, v_n \bigr] = (u_{-\frac{1}{2}}v)_n \qquad \mbox{for $n \in \frac{1}{2}
\mathbb{Z}$},
\end{equation}
and 
\begin{equation}
\bigl[u_{-\frac{1}{2}}, v_{-\frac{1}{2}} \bigr] = (u_{-\frac{1}{2}}v)_{-\frac{1}{2}};
\end{equation}
thus the operators $u_{-\frac{1}{2}}$ form a Lie superalgebra. 

{}From Taylor's Theorem for formal calculus (cf. \cite{FHL}) and the 
$L(-1)$- and $G(-\frac{1}{2})$-derivative properties 
(\ref{L(-1)-derivative}) and (\ref{G(-1/2)-derivative}), we have
\begin{eqnarray}
Y(e^{x_0 L(-1)}v, (x,\varphi)) &=& e^{x_0 \frac{\partial}{\partial x}} 
Y(v,(x,\varphi)) \label{exponential L(-1) property} \\
&=& Y(v,(x_0 + x, \varphi)) \\
Y(e^{\varphi_0 G(-\frac{1}{2})}v, (x,\varphi)) &=& e^{\varphi_0
\left(\frac{\partial}{\partial \varphi} + \varphi
\frac{\partial}{\partial x} \right)} Y(v,(x,\varphi)) \\
&=& Y(v,(x + \varphi_0 \varphi, \varphi_0 + \varphi)), 
\end{eqnarray}
and thus
\begin{eqnarray}
Y(e^{x_0 L(-1) + \varphi_0 G(-\frac{1}{2})}v, (x,\varphi)) &=& e^{x_0
\frac{\partial}{\partial x} + \varphi_0 \left(
\frac{\partial}{\partial \varphi} + \varphi \frac{\partial}{\partial 
x} \right)} Y(v,(x,\varphi)) \label{exponential L(-1) and G(-1/2) property}\\
&=& Y(v,(x_0 + x + \varphi_0 \varphi, \varphi_0 + \varphi)).
\end{eqnarray}

Taking $\mbox{Res}_{x_1}$ of the supercommutator formula (\ref{bracket 
relation for a vosa}) with $u = \tau$ or with $u = \omega =
\frac{1}{2}G(-\frac{1}{2}) \tau$, we have  
\begin{eqnarray}
\Bigl[ L(-1), Y(v,(x,\varphi)) \Bigr] \! &=& \!   Y(L(-1)v,(x,\varphi)) \\
\biggl[ G(-\frac{1}{2}), Y(v,(x,\varphi)) \biggr] \! &=& \! 
Y(G(-\frac{1}{2})v,(x,\varphi)) - 2\varphi Y(L(-1)v,(x,\varphi)) 
\label{G(-1/2) bracket formula} \\
\Bigl[ L(0), Y(v,(x,\varphi)) \Bigr]  \! &=& \! Y(L(0)v,(x,\varphi)) + \frac{\varphi}{2}
Y(G(-\frac{1}{2})v,(x,\varphi)) \\
& & + \; xY(L(-1)v,(x,\varphi))\nonumber\\
\Bigl[ G(\frac{1}{2}), Y(v,(x,\varphi)) \Bigr] \!  &=& \! Y(G(\frac{1}{2})v,(x,\varphi))
- 2\varphi Y(L(0)v,(x,\varphi))\\
& & + \; xY(G(-\frac{1}{2})v,(x,\varphi)) - 2x\varphi Y(L(-1)v,(x,\varphi)) 
\nonumber \\
\Bigl[ L(1), Y(v,(x,\varphi)) \Bigr]  \! &=& \! Y(L(1)v,(x,\varphi)) + \varphi
Y(G(\frac{1}{2})v,(x,\varphi)\\
& & + \; 2x Y(L(0)v,(x,\varphi)) + x\varphi Y(G(-\frac{1}{2})v,(x,\varphi)) \nonumber\\
& & +\; x^2 Y(L(-1)v,(x,\varphi)) .\nonumber
\end{eqnarray}

Making repeated use of properties (\ref{exponential L(-1) property}) 
- (\ref{G(-1/2) bracket formula}), we have
\begin{eqnarray}
e^{x_0 L(-1)} Y(v,(x,\varphi)) e^{-x_0 L(-1)} &=& Y(e^{x_0 L(-1)}v,
(x,\varphi)) \\
&=& Y(v,(x + x_0, \varphi)) \\
e^{\varphi_0 G(-\frac{1}{2})} Y(v,(x,\varphi)) e^{-\varphi_0
G(-\frac{1}{2})} &=& Y(e^{\varphi_0 G(-\frac{1}{2}) - 2\varphi_0
\varphi L(-1)}v,(x,\varphi)) \\
&=& Y(v,(x + \varphi \varphi_0, \varphi + \varphi_0)), 
\end{eqnarray}
and thus
\begin{eqnarray}
& & \hspace{-2.2in} e^{x_0 L(-1) + \varphi_0 G(-\frac{1}{2})} Y(v,(x,\varphi)) e^{- x_0
L(-1) - \varphi_0 G(-\frac{1}{2})} =  \nonumber\\
\hspace{1.9in} &=& Y(e^{x_0 L(-1) + \varphi_0
G(-\frac{1}{2}) - 2\varphi_0 \varphi L(-1)}v,(x,\varphi)) \hspace{.4in}\\
&=& Y(v,(x +x_0  + \varphi \varphi_0, \varphi + \varphi_0)). 
\end{eqnarray} 

{}From the creation property and (\ref{exponential L(-1) and G(-1/2) 
property}), we have
\begin{equation}\label{for skew-symmetry}
e^{x L(-1) + \varphi G(-\frac{1}{2})} v \; = \; Y(v,(x,\varphi))
\mathbf{1} .
\end{equation}

Note that the right-hand side of the Jacobi identity is
invariant under
\[(u,v,x_0,x_1,x_2,\varphi_1,\varphi_2) \longleftrightarrow ((-1)^{\eta(v)
\eta(u)} v,u,- x_0,x_2,x_1,\varphi_2,\varphi_1) .\]
Thus the left-hand side of the Jacobi identity must be symmetric with
respect to this also.  Then using (\ref{for skew-symmetry}), we have
the following {\it skew-supersymmetry} property: for $u,v$
of homogeneous sign in $V$
\begin{equation}
e^{x L(-1) + \varphi G(-\frac{1}{2})} Y(v,(-x,-\varphi))u \; = \;
(-1)^{\eta(v) \eta(u)} Y(u,(x,\varphi))v.
\end{equation}

Let $(V_1, Y_1(\cdot,(x,\varphi)),\mbox{\bf 1}_1,\tau_1)$ and $(V_2,
Y_2(\cdot,(x,\varphi)),\mbox{\bf 1}_2,\tau_2)$ be two vertex operator
superalgebras over $\bigwedge_*$.  A {\it homomorphism of vertex
operator superalgebras with odd formal variables} is a doubly graded
$\bigwedge_*$-module homomorphism $\gamma : V_1 \longrightarrow V_2 \;$
(i.e., $\gamma : (V_1)_{(n)}^i \longrightarrow  (V_2)_{(n)}^i$ for $n \in 
\frac{1}{2} \mathbb{Z}$, and $i \in \mathbb{Z}_2$) such that 
\[\gamma (Y_1(u,(x,\varphi))v) = Y_2(\gamma(u),(x,\varphi))\gamma(v)
\quad \mbox{for} \quad u,v \in V_1 ,\]
$\gamma(\mbox{\bf 1}_1) = \mbox{\bf 1}_2$, and $\gamma(\tau_1) =
\tau_2$.  

\begin{rema}\label{change of square root} {\em In $N=1$ superconformal field
theory there is an inherent choice of square root being made in the choice
of superderivation $D$ satisfying $D^2 = \frac{\partial}{\partial z}$ and
correspondingly in the choice of Neveu-Schwarz algebra representative
element $G(-\frac{1}{2})$ satisfying $G(-\frac{1}{2})^2 = L(-1)$.  In each
case, if $D$ satisfies $D^2 = \frac{\partial}{\partial z}$, then so does
$-D$, and if $G(-\frac{1}{2})$ satisfies $G(-\frac{1}{2})^2 = L(-1)$, then
so does $-G(-\frac{1}{2})$.  The transformation $D \leftrightarrow -D$ 
actually corresponds to the transformation $\varphi \leftrightarrow 
-\varphi$, and algebraically to the transformation $G(-\frac{1}{2}) 
\leftrightarrow -G(-\frac{1}{2})$.  This symmetry gives an isomorphism of 
vertex operator superalgebras with odd formal variables given by
$(V,Y(\cdot,(x,\varphi), \mathbf{1}, \tau) \leftrightarrow (V,Y(\cdot,
(x,-\varphi), \mathbf{1}, -\tau)$.}
\end{rema}

\section{Vertex operator superalgebras over $\bigwedge_*$ and without 
odd formal variables}  

\begin{defn} {\em A} ($N = 1$ Neveu-Schwarz) vertex operator
superalgebra over $\bigwedge_*$  and without odd variables {\em is a
$\frac{1}{2} \mathbb{Z}$-graded (by weight) $\bigwedge_*$-module
which is also $\mathbb{Z}_2$-graded (by sign)
\[V = \coprod_{n \in \frac{1}{2} \mathbb{Z}} V_{(n)} = \coprod_{n \in
\frac{1}{2} \mathbb{Z}} V_{(n)}^0 \oplus \coprod_{n \in \frac{1}{2}
\mathbb{Z}} V_{(n)}^1 = V^0 \oplus V^1 \]  
such that
\[\dim V_{(n)} < \infty \quad \mbox{for} \quad n \in \frac{1}{2}
\mathbb{Z} ,\]
\[V_{(n)} = 0 \quad \mbox{for $n$ sufficiently small} ,\]
equipped with a linear map $V \otimes V \longrightarrow V[[x,x^{-1}]]$, or
equivalently, 
\begin{eqnarray*}
V & \longrightarrow & (\mbox{End} \; V)[[x,x^{-1}]]\\
v & \mapsto & Y(v,x) = \sum_{n \in \mathbb{Z}} v_n x^{-n-1}  
\end{eqnarray*}
where $v_n \in (\mbox{End} \; V)^{\eta(v)}$ for $v$ of homogeneous
sign in $V$, $x$ is an even formal variable, and $Y(v,x)$ denotes the}
vertex operator associated with $v$, {\em and equipped also with two
distinguished homogeneous vectors $\mbox{\bf 1} \in V_{(0)}^0$ (the
{\em vacuum}) and $\tau \in V_{(\frac{3}{2})}^1$ (the {\em Neveu-Schwarz
element}).  The following conditions are assumed for $u,v \in V$: 
\[u_n v = 0 \quad \mbox{for $n \in \mathbb{Z}$ sufficiently large;} \]
\[Y(\mbox{\bf 1}, x) = 1 \quad \mbox{(1 on the right being
the identity operator);} \]
the} creation property {\em holds:
\[Y(v,x) \mbox{\bf 1} \in V[[x]] \qquad \mbox{and} \qquad \lim_{x
\rightarrow 0} Y(v,x) \mbox{\bf 1} = v ; \]
the} Jacobi identity {\em holds: 
\[x_0^{-1} \delta \biggl( \frac{x_1 - x_2}{x_0}\biggr) Y(u,x_1)Y(v,x_2) \hspace{2.8in} \]
\[- (-1)^{\eta(u)\eta(v)} x_0^{-1} \delta \biggl( \frac{ - x_2 +
x_1}{x_0} \biggr)Y(v,x_2)Y(u,x_1) \]
\[ \hspace{2.4in} = x_2^{-1} \delta \biggl( \frac{x_1 - x_0}{x_2} \biggr)
Y(Y(u,x_0)v,x_2) ,\]
for $u,v$ of homogeneous sign in $V$; the Neveu-Schwarz algebra relations hold:
\begin{eqnarray*}
\left[L(m),L(n) \right] \! &=& \! (m - n)L(m + n) + \frac{1}{12} (m^3 - m)
\delta_{m + n , 0} (\mbox{rank} \; V) , \\
\biggl[ G(m + \frac{1}{2}),L(n) \biggr] \! &=& \! (m - \frac{n - 1}{2} ) G(m
+ n + \frac{1}{2}) ,\\
\biggl[ G (m + \frac{1}{2} ) , G(n - \frac{1}{2} ) \biggr] \! &=& \! 2L(m +
n) + \frac{1}{3} (m^2 + m) \delta_{m + n , 0} (\mbox{rank} \; V) ,  
\end{eqnarray*} 
for $m,n \in \mathbb{Z}$, where 
\[G(n + \frac{1}{2}) = \tau_{n + 1} \qquad \mbox{for} \; n \in
\mathbb{Z}, \mbox{i.e.,} \quad Y(\tau,x) = \sum_{n \in \mathbb{Z}} G (n +
\frac{1}{2}) x^{- n - \frac{1}{2} - \frac{3}{2}}, \]  
and $\mbox{rank} \; V \in \mathbb{C}$;
\[L(0)v = nv \quad \mbox{for} \quad n \in \frac{1}{2} \mathbb{Z} \quad
\mbox{and} \quad v \in V_{(n)} ; \]
\begin{equation}\label{L derivative} 
\frac{\partial}{\partial x} Y(v, x) = Y(L(-1)v,x) .
\end{equation}} 
\end{defn}
The vertex operator superalgebra just defined is denoted by 
\[(V,Y(\cdot,x),\mbox{\bf 1},\tau) . \]

Some consequences of the definition are
\begin{equation}\label{L bracket}
[L(- 1), Y(v,x) ] = Y(L(-1)v,x) ;
\end{equation}
\begin{equation}\label{G bracket}
[G(- \frac{1}{2}), Y(v,x) ] = Y(G(- \frac{1}{2})v,x) ;
\end{equation}
\begin{equation}\label{killing of vacuum}
L(n) \mbox{\bf 1} = 0 , \; \mbox{and} \; G(n + \frac{1}{2})
\mbox{\bf 1} = 0 , \; \mbox{for} \; n \geq -1 .
\end{equation}

Note that our definition of vertex operator superalgebra over 
$\bigwedge_*$ (without formal variables) is an obvious extension
of the usual notion of vertex operator algebra \cite{Bo}, 
\cite{FLM}, \cite{FHL} and vertex operator superalgebra (cf. 
\cite{T}, \cite{G}, \cite{FFR}, \cite{DL}, and \cite{KW}). For
instance, any of the examples of ``$N=1$ NS-type SVOAs" can be 
extended in the obvious way to be a $\bigwedge_*$-module 
instead of just a vector space over $\mathbb{C}$, thus giving an 
$N=1$ Neveu-Schwarz vertex operator over $\bigwedge_*$ as defined 
above.   As we shall see in the next section, the categories of 
$N=1$ Neveu-Schwarz vertex operator superalgebras over 
$\bigwedge_*$ with and without odd formal variables,  
respectively, are isomorphic. This isomorphism can be used to give  
examples of $N=1$ Neveu-Schwarz vertex operator algebras over 
$\bigwedge_*$ with odd formal variables {}from the examples of
$N=1$ NS-type SVOAs.

Let $(V_1, Y_1(\cdot,x),\mbox{\bf 1}_1,\tau_1)$ and $(V_2,
Y_2(\cdot,x),\mbox{\bf 1}_2,\tau_2)$ be two vertex operator 
superalgebras over $\bigwedge_*$.  A {\it homomorphism of vertex
operator superalgebras without odd formal variables} is a
doubly graded $\bigwedge_*$-module homomorphism $\gamma : V_1
\longrightarrow V_2$ such that 
\[\gamma (Y_1(u,x)v) = Y_2(\gamma(u),x)\gamma(v) \quad \mbox{for}
\quad u,v \in V_1 ,\] 
$\gamma(\mbox{\bf 1}_1) = \mbox{\bf 1}_2$, and $\gamma(\tau_1) =
\tau_2$.

\section{The isomorphism between the category of vertex operator
superalgebras with odd formal variables and the category of vertex
operator superalgebras without odd formal variables}

\begin{prop}\label{get a vosa without}
Let $(V,Y(\cdot,(x,\varphi)), \mbox{\bf 1}, \tau)$ be a vertex
operator superalgebra with odd formal variables.  Then
$(V,Y(\cdot,(x,0)), \mbox{\bf 1}, \tau)$ is a vertex operator
superalgebra without odd formal variables. 
\end{prop}

{\it Proof:} \hspace{.2cm}  
It is trivial that $(V,Y(\cdot,(x,0)),
\mbox{\bf 1}, \tau)$ satisfies all the axioms for a vertex operator
superalgebra without odd formal variables except for the
$L(-1)$-derivative property (\ref{L derivative}).  But this follows as
a consequence of the definition of a vertex operator superalgebra with
odd formal variables (\ref{L(-1)-derivative}). $\pfbox$ \\

Let $(V,Y(\cdot,x), \mathbf{1}, \tau)$ be a vertex operator
superalgebra without odd formal variables.  Define
\[\tilde{Y}(v, (x,\varphi)) = Y(v,x) + \varphi Y(G(-\frac{1}{2}) v,x) .\]

\begin{prop}\label{get a vosa with}
$(V,\tilde{Y}(\cdot,(x,\varphi)), \mbox{\bf 1}, \tau)$ is a vertex
operator superalgebra with odd formal variables. 
\end{prop}

{\it Proof:} \hspace{.2cm} Axioms (\ref{vosa2}), (\ref{positive
energy}), (\ref{truncation}), and (\ref{grading for vosa with}) are
obvious.  Axiom (\ref{vacuum identity}) holds since by consequence
(\ref{killing of vacuum}) of the definition of $(V,Y(\cdot,x),
\mathbf{1}, \tau)$, we have 
\[\tilde{Y}(\mathbf{1}, (x.\varphi)) = Y(\mathbf{1},x) + \varphi
Y(G(-\frac{1}{2}) \mathbf{1}, x) = Y(\mathbf{1},x) + \varphi Y(0, x) = 1 .\] 
The creation property holds since by the creation property for
$(V,Y(\cdot,x), \mathbf{1}, \tau)$ 
\[\tilde{Y}(v,(x,\varphi)) \mathbf{1} = Y(v,x)\mathbf{1} + \varphi
Y(G(-\frac{1}{2}) v, x) \mathbf{1} \in V[[x]] [\varphi] ,\]
and 
\begin{eqnarray*}
\lim_{(x,\varphi) \rightarrow 0} \tilde{Y}(v,(x,\varphi)) \mathbf{1}
&=& \lim_{(x,\varphi) \rightarrow 0} Y(v,x)\mathbf{1} +
\lim_{(x,\varphi) \rightarrow 0} \varphi  Y(G(-\frac{1}{2}) v, x) \mathbf{1}
\\
&=&  \lim_{x \rightarrow 0} Y(v,x)\mathbf{1} = v .
\end{eqnarray*}

To prove the Jacobi identity for $(V,\tilde{Y}(\cdot,(x,\varphi)),
\mbox{\bf 1}, \tau)$, we use the Jacobi identity for $(V,Y(\cdot,x),
\mathbf{1}, \tau)$ and the $L(-1)$- and $G(-\frac{1}{2})$-bracket properties
for $(V,Y(\cdot,x), \mathbf{1}, \tau)$ given by (\ref{L bracket}) and 
(\ref{G bracket}).  
\begin{eqnarray*}
& & \hspace{-.4in} x_0^{-1} \delta \biggl( \frac{x_1 - x_2 - \varphi_1 \varphi_2}{x_0}
\biggr) \tilde{Y}(u,(x_1, \varphi_1)) \tilde{Y}(v,(x_2, \varphi_2)) \\
& &  - (-1)^{\eta(u)\eta(v)} x_0^{-1} \delta \biggl(
\frac{x_2 - x_1 + \varphi_1 \varphi_2}{-x_0} \biggr) \tilde{Y}(v,(x_2,
\varphi_2)) \tilde{Y}(u,(x_1, \varphi_1)) \\
&=& x_0^{-1} \delta \biggl( \frac{x_1 - x_2 - \varphi_1 \varphi_2}{x_0}
\biggr) \biggl(Y(u,x_1) + \varphi_1 Y(G(-\frac{1}{2})u,x_1) \biggr) \\
& & \hspace{2.7in} \biggl( Y(v,x_2) + \varphi_2 Y(G(-\frac{1}{2})v,x_2)
\biggr)\\ 
& & - (-1)^{\eta(u)\eta(v)} x_0^{-1} \delta \biggl( \frac{x_2 - x_1 + 
\varphi_1 \varphi_2}{-x_0} \biggr)  \biggl( Y(v,x_2) + \varphi_2
Y(G(-\frac{1}{2})v,x_2) \biggr) \\
& & \hspace{2.7in} \biggl(Y(u,x_1) + \varphi_1 Y(G(-\frac{1}{2})u,x_1)
\biggr) \\ 
&=& x_0^{-1} \delta \biggl( \frac{x_1 - x_2}{x_0}\biggr) Y(u,x_1)
Y(v,x_2) \\
& & \hspace{1.8in} - (-1)^{\eta(u)\eta(v)} x_0^{-1} \delta \biggl( \frac{x_2 -
x_1}{-x_0} \biggr) Y(v,x_2) Y(u,x_1) \\
& & + \; x_0^{-1} \delta \biggl( \frac{x_1 - x_2}{x_0}\biggr) \varphi_1 
Y(G(-\frac{1}{2})u,x_1)Y(v,x_2) \\
& & \hspace{1.1in} - (-1)^{\eta(u)\eta(v)} x_0^{-1} \delta \biggl(
\frac{x_2 - x_1}{-x_0} \biggr) Y(v,x_2) \varphi_1
Y(G(-\frac{1}{2})u,x_1) \\ 
& & + \; x_0^{-1} \delta \biggl( \frac{x_1 - x_2}{x_0}\biggr)  
Y(u,x_1)\varphi_2 Y(G(-\frac{1}{2})v,x_2) \\
& & \hspace{1.1in} - (-1)^{\eta(u)\eta(v)} x_0^{-1} \delta \biggl(
\frac{x_2 - x_1}{-x_0} \biggr)\varphi_2 Y(G(-\frac{1}{2})v,x_2)
Y(u,x_1) \\ 
& & + \; x_0^{-1} \delta \biggl( \frac{x_1 - x_2}{x_0}\biggr) \varphi_1
Y(G(-\frac{1}{2})u,x_1)\varphi_2 Y(G(-\frac{1}{2})v,x_2)\\
& & \hspace{.4in} - (-1)^{\eta(u)\eta(v)} x_0^{-1} \delta \biggl(
\frac{x_2 - x_1}{-x_0} \biggr)\varphi_2 Y(G(-\frac{1}{2})v,x_2)
\varphi_1 Y(G(-\frac{1}{2})u,x_1) \\ 
& & - \varphi_1 \varphi_2 x_0^{-2} \delta' \biggl( \frac{x_1 -
x_2}{x_0}\biggr) Y(u,x_1) Y(v,x_2) \\
& & \hspace{1.3in} + \; \varphi_1 \varphi_2 (-1)^{\eta(u)\eta(v)} x_0^{-2}
\delta' \biggl( \frac{x_2 - x_1}{-x_0} \biggr) Y(v,x_2) Y(u,x_1) \\ 
&=& x_2^{-1} \delta \biggl( \frac{x_1 - x_0}{x_2}\biggr)
Y(Y(u,x_0)v,x_2) \\
& & + \; \varphi_1 x_0^{-1} \delta \biggl( \frac{x_1 -
x_2}{x_0}\biggr) Y(G(-\frac{1}{2})u,x_1)Y(v,x_2) \\
& &  - \; \varphi_1 (-1)^{\eta(u)\eta(v) + \eta(v)} x_0^{-1}
\delta \biggl( \frac{x_2 - x_1}{-x_0} \biggr) Y(v,x_2) 
Y(G(-\frac{1}{2})u,x_1) \\ 
& & + \; \varphi_2 (-1)^{\eta(u)} x_0^{-1} \delta \biggl( \frac{x_1 -
x_2}{x_0}\biggr) Y(u,x_1)Y(G(-\frac{1}{2})v,x_2) \\
& &  - \; \varphi_2 (-1)^{\eta(u)\eta(v)} x_0^{-1} \delta
\biggl( \frac{x_2 - x_1}{-x_0} \biggr) Y(G(-\frac{1}{2})v,x_2)
Y(u,x_1) \\ 
& &  + \; \varphi_1 \varphi_2 (-1)^{\eta(u) + 1} x_0^{-1} \delta \biggl(
\frac{x_1 - x_2}{x_0} \biggr) Y(G(-\frac{1}{2})u,
x_1)Y(G(-\frac{1}{2})v,x_2) \\
& & - \; \varphi_1 \varphi_2 (-1)^{\eta(u)\eta(v) +
\eta(v)} x_0^{-1} \delta \biggl( \frac{x_2 - x_1}{-x_0}
\biggr)Y(G(-\frac{1}{2})v,x_2) Y(G(-\frac{1}{2})u,x_1) \\ 
& & - \; \varphi_1 \varphi_2 \frac{\partial}{\partial x_1} \left(
x_0^{-2} \delta \biggl( \frac{x_1 - x_2}{x_0}\biggr) \right) Y(u,x_1)
Y(v,x_2) \\ 
& & + \; \varphi_1 \varphi_2 (-1)^{\eta(u)\eta(v)}
\frac{\partial}{\partial x_1} \left( x_0^{-2} \delta \biggl( \frac{x_2
- x_1}{-x_0} \biggr) \right) Y(v,x_2) Y(u,x_1) \\  
&=& x_2^{-1} \delta \biggl( \frac{x_1 - x_0}{x_2}\biggr)
Y(Y(u,x_0)v,x_2) \\
& & + \; \varphi_1 x_2^{-1} \delta \biggl( \frac{x_1 -
x_0}{x_2}\biggr) Y(Y(G(-\frac{1}{2})u,x_0)v,x_2) \\
& & + \; \varphi_2 (-1)^{\eta(u)} x_2^{-1} \delta \biggl( \frac{x_1 -
x_0}{x_2}\biggr) Y(Y(u,x_0)G(-\frac{1}{2})v,x_2) \\
& & + \; \varphi_1 \varphi_2 (-1)^{\eta(u) + 1} x_2^{-1} \delta \biggl(
\frac{x_1 - x_0}{x_2} \biggr) Y(Y(G(-\frac{1}{2})u,
x_0)G(-\frac{1}{2})v,x_2) \\  
& & - \; \varphi_1 \varphi_2 \frac{\partial}{\partial x_1} \left(
x_0^{-2} \delta \biggl( \frac{x_1 - x_2}{x_0}\biggr) Y(u,x_1) Y(v,x_2)
\right. \\
& &  \hspace{1.4in} \left. - (-1)^{\eta(u)\eta(v)} x_0^{-2} \delta \biggl(
\frac{x_2 - x_1}{-x_0} \biggr) Y(v,x_2) Y(u,x_1) \right) \\  
& & + \; \varphi_1 \varphi_2 \left(x_0^{-2} \delta \biggl( \frac{x_1 -
x_2}{x_0}\biggr) \frac{\partial}{\partial x_1} Y(u,x_1) Y(v,x_2)
\right. \\
& & \left. \hspace{1.4in} - (-1)^{\eta(u)\eta(v)} x_0^{-2} \delta \biggl(
\frac{x_2 - x_1}{-x_0} \biggr) Y(v,x_2) \frac{\partial}{\partial x_1}
Y(u,x_1) \right) \\ 
&=& x_2^{-1} \delta \biggl( \frac{x_1 - x_0}{x_2}\biggr) \left(
Y(Y(u,x_0)v,x_2) + \varphi_1 Y(Y(G(-\frac{1}{2})u,x_0)v,x_2)
\right. \\ 
& & \hspace{.6in} + \; \varphi_2 (-1)^{\eta(u)} Y(Y(u,x_0)
G(-\frac{1}{2})v, x_2) \\
& & \hspace{1.5in} \left. + \; \varphi_1 \varphi_2 (-1)^{\eta(u) + 1}
Y(Y(G(-\frac{1}{2})u, x_0) G(-\frac{1}{2})v,x_2) \right)\\  
& & - \; \varphi_1 \varphi_2 \frac{\partial}{\partial x_1} \left(
x_2^{-1} \delta \biggl( \frac{x_1 - x_0}{x_2}\biggr) Y(Y(u,x_0)v,x_2)
\right) \\
& & + \; \varphi_1 \varphi_2 \left(x_0^{-2} \delta \biggl( \frac{x_1 -
x_2}{x_0}\biggr) Y(L(-1)u,x_1) Y(v,x_2)
\right. \\
& & \hspace{1.3in} \left. - (-1)^{\eta(u)\eta(v)} x_0^{-2} \delta \biggl(
\frac{x_2 - x_1}{-x_0} \biggr) Y(v,x_2) Y(L(-1)u,x_1) \right) \\
&=& x_2^{-1} \delta \biggl( \frac{x_1 - x_0}{x_2}\biggr) \Biggl(
Y(Y(u,x_0)v,x_2) + \varphi_1 Y(Y(G(-\frac{1}{2})u,x_0)v,x_2)
\Biggr. \\
& & + \; \varphi_2 \biggl( (-1)^{\eta(u)}
Y(Y(u,x_0)G(-\frac{1}{2})v,x_2) + Y(\biggl[G(-\frac{1}{2}), Y(u,x_0) 
\biggr]v, x_2) \biggr. \\
& & \hspace{3in} \biggl. - Y(Y(G(-\frac{1}{2})u,x_0)v,x_2) \biggr) \\
& &  + \; \varphi_1 \varphi_2 \biggl( (-1)^{\eta(u) +
1} Y(Y(G(-\frac{1}{2})u,x_0)G(-\frac{1}{2})v,x_2) \biggr.\\
& & \hspace{3in} \Biggl. \biggl. + Y(Y(L(-1)u,x_0)v,
x_2) \biggr) \Biggr) \\
& & - \; \varphi_1 \varphi_2 x_2^{-2} \delta '\biggl( \frac{x_1 -
x_0}{x_2}\biggr) Y(Y(u,x_0)v,x_2) \\
&=& x_2^{-1} \delta \biggl( \frac{x_1 - x_0 - \varphi_1
\varphi_2}{x_2}\biggr) \Biggl( Y(Y(u,x_0)v,x_2) + \; \varphi_1
Y(Y(G(-\frac{1}{2})u,x_0)v,x_2) \Biggr. \\
& &+ \; \varphi_2 Y(G(-\frac{1}{2})Y(u,x_0)v, x_2) - \varphi_2
Y(Y(G(-\frac{1}{2})u,x_0)v,x_2) \\
& & + \; \varphi_1 \varphi_2 \biggl( (-1)^{\eta(u) + 1}
Y(Y(G(-\frac{1}{2})u, x_0) G(-\frac{1}{2})v,x_2) 
\biggr. \\ 
& & \hspace{2.2in} \Biggl. + \biggl. \; Y(\biggl[G(-\frac{1}{2}),
Y(G(-\frac{1}{2})u,x_0) \biggr] v,x_2) \biggr) \Biggr) \\  
&=& x_2^{-1} \delta \biggl( \frac{x_1 - x_0 - \varphi_1
\varphi_2}{x_2}\biggr) \left( Y \Biggl(Y(u,x_0)v + \varphi_1
Y(G(-\frac{1}{2})u,x_0)v \Biggr. \right.\\
& & \Biggl.  - \; \varphi_2 Y(G(-\frac{1}{2})u,x_0)v ,x_2 \Biggr)  + \;
\varphi_2  Y(G(-\frac{1}{2})Y(u,x_0)v, x_2) \\
& & \hspace{2.2in} \left. + \; \varphi_1 \varphi_2
Y(G(-\frac{1}{2})Y(G(-\frac{1}{2})u,x_0)v,x_2) \right) \\ 
&=& x_2^{-1} \delta \biggl( \frac{x_1 - x_0 - \varphi_1
\varphi_2}{x_2}\biggr) \Biggl( Y \biggl( \tilde{Y}(u, (x_0,\varphi_1 -
\varphi_2) )v, x_2 \biggr)  \Biggr. \\
& & \hspace{2.1in} \Biggl. + \; \varphi_2 Y(G(-\frac{1}{2})\tilde{Y}(u,
(x_0,\varphi_1 - \varphi_2) )v, x_2) \Biggr)\\
&=& x_2^{-1} \delta \biggl( \frac{x_1 - x_0 - \varphi_1
\varphi_2}{x_2}\biggr) \tilde{Y}( \tilde{Y}(u, (x_0,\varphi_1 -
\varphi_2) )v, (x_2, \varphi_2)) \\
\end{eqnarray*}
which gives the Jacobi identity for $(V,\tilde{Y}(\cdot,(x,\varphi)),
\mathbf{1}, \tau)$.  

For the Neveu-Schwarz element $\tau$, we have
\begin{eqnarray*}
\tilde{Y} (\tau, (x,\varphi)) &=& Y(\tau,x) + \varphi
Y(G(-\frac{1}{2})\tau,x) \\
&=& Y(\tau,x) + \varphi \biggl[G(-\frac{1}{2}), Y(\tau,x) \biggr]\\
&=& \sum_{n \in \mathbb{Z}} G(n + \frac{1}{2}) x^{-n-2} + \varphi \sum_{n
\in \mathbb{Z}} \biggl[G(-\frac{1}{2}), G(n + \frac{1}{2}) \biggr]
x^{-n-2}  \\
&=& \sum_{n \in \mathbb{Z}} G(n + \frac{1}{2}) x^{-n-2} + \varphi \sum_{n
\in \mathbb{Z}} 2 L(n) x^{-n-2}  
\end{eqnarray*}
which gives (\ref{stress tensor}).  Finally, using the
$L(-1)$-derivative property for $(V,Y(\cdot,x),\mathbf{1}, \tau)$, we
have
\begin{eqnarray*}
\biggl(\frac{\partial}{\partial \varphi} + \varphi
\frac{\partial}{\partial x}  \biggr) \tilde{Y} (v, (x ,\varphi)) &=&
\biggl(\frac{\partial}{\partial \varphi} + \varphi 
\frac{\partial}{\partial x}  \biggr) \biggl( Y (v, x) + \varphi Y
(G(-\frac{1}{2})v, x) \biggr) \\
&=& \varphi \frac{\partial}{\partial x}Y (v, x) + Y(G(-\frac{1}{2})v,
x) \\
&=& Y(G(-\frac{1}{2})v,x) + \varphi Y (L(-1)v, x) \\
&=& Y(G(-\frac{1}{2})v,x) + \varphi Y (G(-\frac{1}{2})^2 v, x) \\
&=& \tilde{Y}(G(-\frac{1}{2})v,(x,\varphi))
\end{eqnarray*} 
which gives the $G(-\frac{1}{2})$-derivative property
(\ref{G(-1/2)-derivative}). Thus
$(V,\tilde{Y}(\cdot,(x,\varphi)),\mathbf{1}, \tau)$ is a vertex
operator superalgebra with odd formal variables.   $\pfbox$ \\

Let $\mathbf{SV}(\varphi,c,*)$ denote the category of vertex operator
superalgebras over $\bigwedge_*$ with odd formal variables and with
central charge, i.e., rank, $c \in \mathbb{C}$, and let $\mathbf{SV}(c,*)$ 
denote the category of vertex operator superalgebras over $\bigwedge_*$ 
without odd formal variables and with central charge, i.e., rank, $c \in 
\mathbb{C}$.  Let $1_{SV_\varphi}$ and $1_{SV}$ be the identity functors 
on the categories $\mathbf{SV}(\varphi,c,*)$ and $\mathbf{SV}(c,*)$, 
respectively. 

\begin{thm}\label{superalgebras} 
For any $c \in \mathbb{C}$, the two categories $\mathbf{SV}(\varphi,c,*)$
and $\mathbf{SV}(c,*)$ are isomorphic.  That is there exist two
functors $F_0: \mathbf{SV}(\varphi,c,*) \longrightarrow \mathbf{SV}(c,*)$
and $F_\varphi : \mathbf{SV}(c,*) \longrightarrow
\mathbf{SV}(\varphi,c,*)$ such that $F_0 \circ F_\varphi = 1_{SV}$ and
$F_\varphi \circ F_0 = 1_{SV_\varphi}$. 
\end{thm}

{\it Proof:} \hspace{.2cm} We first define $F_0$ by 
\[F_0(V,Y(\cdot,(x,\varphi)), \mathbf{1}, \tau) = (V,Y(\cdot,(x,0)),
\mathbf{1}, \tau), \quad \mbox{and} \quad F_0(\gamma) = \gamma .\]
Proposition \ref{get a vosa without} shows that $F_0$ takes objects in 
$\mathbf{SV}(\varphi,c,*)$ to objects in $\mathbf{SV}(c,*)$.  It is
clear that $F_0$ takes morphisms in $\mathbf{SV}(\varphi,c,*)$ to 
morphisms in $\mathbf{SV}(c,*)$ and that $F_0$ is a functor.  

We next define $F_\varphi$ by 
\[F_\varphi(V,Y(\cdot,x), \mathbf{1}, \tau) = (V,\tilde{Y} (\cdot,
(x,\varphi)), \mathbf{1}, \tau), \quad \mbox{and} \quad 
F_\varphi(\gamma) = \gamma\] 
where $\tilde{Y}(v,(x,\varphi)) = Y(v,x) + \varphi Y(G(-
\frac{1}{2})v, x)$.  Proposition \ref{get a vosa with} shows that
$F_\varphi$ takes objects in $\mathbf{SV}(c,*)$ to objects in
$\mathbf{SV} (\varphi,c,*)$.  Let  
\[\gamma: (V_1,Y_1(\cdot,x), \mathbf{1}_1, \tau_1) \longrightarrow (V_2,
Y_2(\cdot,x), \mathbf{1}_2, \tau_2)\]
be a homomorphism of vertex operator superalgebras without odd formal
variables.  Then $\gamma(\tau_1) = \tau_2$.  Denoting $(\tau_1)_0 =
G_1(-\frac{1}{2})$ and $(\tau_2)_0 = G_2(-\frac{1}{2})$, we have
$\gamma(G_1(-\frac{1}{2})u) = G_2(-\frac{1}{2})\gamma(u)$.  Thus  
\begin{eqnarray*}
\gamma(\tilde{Y}_1(u,(x,\varphi))v) &=& \gamma(Y_1(u,x)v + \varphi
Y_1(G_1(-\frac{1}{2})u,x)v ) \\
&=& \gamma(Y_1(u,x)v)  + \varphi \gamma(Y_1(G_1(-\frac{1}{2})u,x)v ) \\
&=& Y_2(\gamma(u), x)\gamma(v) + \varphi
Y_1(\gamma(G_1(-\frac{1}{2})u),x)\gamma(v) \\
&=& Y_2(\gamma(u), x)\gamma(v)+ \varphi
Y_1(G_2(-\frac{1}{2})\gamma(u),x)\gamma(v) \\
&=& \tilde{Y}_2(\gamma(u), (x,\varphi))\gamma(v)
\end{eqnarray*} 
which shows that $F_\varphi(\gamma) = \gamma$ is a homomorphism of
vertex operator superalgebras with odd formal variables, i.e.,
$F_\varphi$ takes morphisms in $\mathbf{SV}(c,*)$ to morphisms in
$\mathbf{SV}(\varphi,c,*)$.  It is clear that $F_\varphi$ is a 
functor. 

The fact that $F_0 \circ F_\varphi = 1_{SV}$ and $F_\varphi \circ F_0 =
1_{SV_\varphi}$ on morphisms is trivial.  The fact that $F_0 \circ
F_\varphi = 1_{SV}$ on objects is given by 
\[Y(v,(x,\varphi)) = Y(v,(x,0)) + \varphi Y(G(-\frac{1}{2})v,(x,0)) ,
\] 
and the fact that $F_\varphi \circ F_0 = 1_{SV_\varphi}$ on objects is 
given by 
\[Y(v,x) = \left. Y(v,x) + \varphi Y(G(-\frac{1}{2})v,x)
\right|_{\varphi = 0} .\] $\pfbox$ \\

\section{Weak supercommutativity and weak associativity for vertex
operator superalgebras with odd formal variables}

In this section we show that the properties which we call ``weak'' 
supercommutativity and ``weak'' associativity for vertex operator 
superalgebras without odd formal variables, as formulated and studied 
for instance in \cite{DL}, \cite{L1} and \cite{L2}, 
have the expected analogues when we add the odd formal variables.  We 
refer to these properties as ``weak'' supercommutativity and ``weak''
associativity because in Section 9, we will prove slightly
stronger statements about the nature of certain rational functions
associated with products and iterates of vertex operators.  We also
note that we could prove Propositions \ref{weak commutativity},
\ref{weak associativity} and \ref{equals Jacobi1} by using weak
commutativity, weak associativity and their equivalence with the
Jacobi identity, respectively, for a vertex operator superalgebra
without odd formal variables and then by using Theorem
\ref{superalgebras}.  However, we choose to prove these properties
directly using the definition of vertex operator superalgebra with odd
formal variables.

In this section we are following and extending the corresponding 
results and arguments of \cite{L1} and \cite{L2}.

Let $\Z$ denote the positive integers.  

\begin{prop}\label{weak commutativity}
{\bf (weak supercommutativity)} Let $(V,Y(\cdot,(x, \varphi)),
\mathbf{1}, \tau)$ be a vertex operator algebra with odd formal
variables and $u,v \in V$ with homogeneous sign.  Then there exists $k
\in \Z$ such that 
\[(x_1 - x_2 - \varphi_1 \varphi_2)^k Y(u,(x_1,\varphi_1))
Y(v,(x_2,\varphi_2)) =  \hspace{1.9in} \]
\[\hspace{1in} (-1)^{\eta(u) \eta(v)} (x_1 - x_2 - \varphi_1
\varphi_2)^k Y(v,(x_2,\varphi_2)) Y(u,(x_1,\varphi_1)) .\]
Furthermore this weak supercommutativity follows {}from the truncation
condition (\ref{truncation}) and the Jacobi identity. 
\end{prop}  

{\it Proof:} \hspace{.2cm} Let $m \in \Z$.  Taking $\mbox{Res}_{x_0}
x_0^{m}$ of the Jacobi identity, we have
\begin{eqnarray*}
& & \hspace{-.4in} (x_1 - x_2 - \varphi_1 \varphi_2)^m \left[ Y(u,(x_1,\varphi_1)),
Y(v,(x_2,\varphi_2)) \right] = \\
&=& (x_1 - x_2 - \varphi_1 \varphi_2)^m
\Bigl(Y(u,(x_1,\varphi_1))Y(v,(x_2,\varphi_2)) \Bigr. \\
& & \hspace{1.9in} \Bigl. - (-1)^{\eta(u) \eta(v)}Y(v,(x_2,\varphi_2))
Y(u,(x_1,\varphi_1)) \Bigr)\\ 
&=& \mbox{Res}_{x_0} x_2^{-1} \delta \biggl( \frac{x_1 - x_0 -
\varphi_1 \varphi_2 }{x_2}\biggr) x_0^{m} Y(Y(u,x_0,\varphi_1 -
\varphi_2))v, (x_2, \varphi_2)) \\
&=& \mbox{Res}_{x_0} x_1^{-1}  \delta \biggl( \frac{x_2 + x_0 + \varphi_1
\varphi_2 }{x_1}\biggr) x_0^{m} Y(Y(u,x_0,\varphi_1 - \varphi_2))v,
(x_2, \varphi_2)) \\
&=& \sum_{n \in \mathbb{N}} \frac{1}{n !} \Biggl(
\biggl(\frac{\partial}{\partial x_2} \biggr)^n x_1^{-1} \delta
\biggl(\frac{x_2}{x_1}\biggr)\Biggr) Y((u_{n+m} + (\varphi_1 -
\varphi_2) u_{n+m - \frac{1}{2}})v, (x_2, \varphi_2)) \\
& & + \; \varphi_1 \varphi_2 \sum_{n \in \mathbb{N}}
\frac{1}{n !} \Biggl( \biggl(\frac{\partial}{\partial x_2} \biggr)^{n +
1} x_1^{-1} \delta \biggl(\frac{x_2}{x_1}\biggr)\Biggr) \\
& & \hspace{1.8in} Y((u_{n+m - 1} + (\varphi_1 - \varphi_2)u_{n+m -
\frac{3}{2}})v, (x_2, \varphi_2)) . 
\end{eqnarray*}
Let $k \in \Z$ be such that $u_l v = 0$ for all $l \in \frac{1}{2}
\Z$, $l \geq k - \frac{3}{2}$.  Setting $m = k$, we obtain weak
supercommutativity.  $\pfbox$ \\

\begin{prop}\label{weak associativity}
{\bf (weak associativity)} Let $(V,Y(\cdot,(x, \varphi)), \mathbf{1},
\tau)$ be a vertex operator algebra with odd formal variables and $u,v
\in V$ with homogeneous sign.  Then there exists $k \in \Z$ such that
for any $w \in V$ 
\[(x_0 + x_2 + \varphi_1 \varphi_2)^k Y(Y(u,(x_0,\varphi_1 -
\varphi_2)v, (x_2,\varphi_2))w =  \hspace{1.5in} \]  
\[\hspace{1.1in} (x_0 + x_2 + \varphi_1 \varphi_2)^k Y(u,(x_0 + x_2 +
\varphi_1 \varphi_2 ,\varphi_1)) Y(v,(x_2,\varphi_2))w  \]
Furthermore this weak associativity follows {}from the truncation
condition (\ref{truncation}) and the Jacobi identity. 
\end{prop}  

{\it Proof:} \hspace{.2cm} Taking $\mbox{Res}_{x_1}$ of the Jacobi
identity, we obtain the following iterate
\begin{eqnarray*}
& & \hspace{-.4in} Y(Y(u,(x_0,\varphi_1 - \varphi_2)v, (x_2,\varphi_2)) =\\
&=& \mbox{Res}_{x_1} x_1^{-1} \delta \biggl( \frac{x_2 + x_0 + \varphi_1
\varphi_2 }{x_1}\biggr) Y(Y(u,(x_0,\varphi_1 - \varphi_2)v,
(x_2,\varphi_2)) \\
&=& \mbox{Res}_{x_1} x_2^{-1} \delta \biggl( \frac{x_1 - x_0 - \varphi_1
\varphi_2 }{x_2}\biggr) Y(Y(u,(x_0,\varphi_1 - \varphi_2)v,
(x_2,\varphi_2)) \\
&=& \mbox{Res}_{x_1} \left( x_0^{-1} \delta \biggl( \frac{x_1 - x_2 -
\varphi_1 \varphi_2}{x_0} \biggr) Y(u,(x_1, \varphi_1))Y(v,(x_2,
\varphi_2)) \right. \\
& & \quad \left. - (-1)^{\eta(u)\eta(v)} x_0^{-1} \delta \biggl(
\frac{x_2 - x_1 + \varphi_1 \varphi_2}{-x_0} \biggr)Y(v,(x_2,
\varphi_2))Y(u,(x_1, \varphi_1)) \right) \\
&=& \mbox{Res}_{x_1} \left( x_1^{-1} \delta \biggl( \frac{x_0 + x_2 +
\varphi_1 \varphi_2}{x_1} \biggr) Y(u,(x_1, \varphi_1))Y(v,(x_2,
\varphi_2)) \right. \\
& & \quad - (-1)^{\eta(u)\eta(v)} Y(v,(x_2, \varphi_2))
\left(x_0^{-1} \delta \biggl(\frac{x_1 - x_2 - \varphi_1
\varphi_2}{x_0} \biggr) \right. \\
& & \hspace{1.7in}\left. \left. - x_2^{-1} \delta \biggl(\frac{x_1 -
x_0 - \varphi_1 \varphi_2}{x_2} \biggr) \right) Y(u,(x_1, \varphi_1))
\right) \\ 
&=& Y(u,(x_0 + x_2 + \varphi_1 \varphi_2 ,\varphi_1))
Y(v,(x_2,\varphi_2)) \\
& & \quad - (-1)^{\eta(u)\eta(v)} Y(v,(x_2, \varphi_2))
\Bigl(Y(u,(x_0 + x_2 + \varphi_1 \varphi_2 ,\varphi_1)) \Bigr. \\
& & \hspace{2.6in} \Bigl. - Y(u,(x_2 + x_0 + \varphi_1 \varphi_2
,\varphi_1)) \Bigr) .\\  
\end{eqnarray*}
For any $w \in V$, let $k \in \Z$ be such that $x^k Y(u,(x,\varphi))w$
involves only positive powers of $x$.  Then 
\begin{eqnarray*}
(x_0 + x_2 + \varphi_1 \varphi_2)^k \Bigl(Y(u,(x_0 + x_2 + \varphi_1
\varphi_2 ,\varphi_1)) & & \\
& & \hspace{-.9in} - \; Y(u,(x_2 + x_0 + \varphi_1 \varphi_2
,\varphi_1)) \Bigr)w   = 0 
\end{eqnarray*}
and weak associativity follows. $\pfbox$ \\

\begin{prop}\label{equals Jacobi1}  
In the presence of the other axioms in the definition of vertex
operator superalgebra with odd formal variables, the Jacobi identity
is equivalent to weak supercommutativity and weak associativity. 
\end{prop}

{\it Proof:} \hspace{.2cm} Propositions \ref{weak commutativity}
and \ref{weak associativity} show that in the presence of the
other axioms for a vertex operator superalgebra with odd formal
variables, the Jacobi identity implies weak supercommutativity and
weak associativity. 

Assume weak supercommutativity and weak associativity hold.  Choose $k
\in \Z$ such that $u_m v = u_m w = 0$ for all $m \in \frac{1}{2} \Z$,
$m \geq k$.  Then
\begin{eqnarray}
& & \hspace{-.4in} x_0^k x_1^k \left(x_0^{-1} \delta \biggl(\frac{x_1 - x_2 - \varphi_1
\varphi_2}{x_0} \biggr)  Y(u,(x_1, \varphi_1))Y(v,(x_2, \varphi_2))w 
\right. \label{assoc1} \\
& & \left. - (-1)^{\eta(u)\eta(v)} x_0^{-1} \delta
\biggl(\frac{x_2 - x_1 + \varphi_1 \varphi_2}{-x_0} \biggr) Y(v,(x_2,
\varphi_2)) Y(u,(x_1, \varphi_1))w \right) \nonumber\\  
&=& x_0^{-1} \delta \biggl(\frac{x_1 - x_2 - \varphi_1 \varphi_2}{x_0}
\biggr) x_1^k (x_1 - x_2 - \varphi_1 \varphi_2)^k \nonumber \\
& &  \hspace{2.7in} Y(u,(x_1, \varphi_1))Y(v,(x_2, \varphi_2))w \nonumber \\
& &  - (-1)^{\eta(u)\eta(v)} x_0^{-1} \delta
\biggl(\frac{x_2 - x_1 + \varphi_1 \varphi_2}{-x_0} \biggr) x_1^k (x_1
- x_2 - \varphi_1 \varphi_2)^k \nonumber\\
& & \hspace{2.7in} Y(v,(x_2, \varphi_2)) Y(u,(x_1, \varphi_1))w \nonumber \\
&=& x_2^{-1} \delta \biggl(\frac{x_1 - x_0 - \varphi_1 \varphi_2}{x_2}
\biggr) \Bigl(x_1^k (x_1 - x_2 - \varphi_1 \varphi_2)^k \Bigr. \nonumber\\
& & \Bigl. \hspace{2.7in} Y(u,(x_1,
\varphi_1))Y(v,(x_2, \varphi_2))w \Bigr) .\nonumber
\end{eqnarray}
By weak supercommutativity, $x_1^k (x_1 - x_2 - \varphi_1 \varphi_2)^k
Y(u,(x_1, \varphi_1))Y(v,(x_2, \varphi_2))w$ involves only nonnegative
powers of $x_1$, and $u_m w = 0$ for $m \geq k$.  Thus in this case,
we can replace $x_1$ by $x_2 + x_0 + \varphi_1 \varphi_2$ or $x_0 +x_2
+\varphi_1 \varphi_2$.  Therefore (\ref{assoc1}) is equal to
\begin{eqnarray*}
\lefteqn{x_2^{-1} \delta \biggl(\frac{x_1 - x_0 - \varphi_1 \varphi_2}{x_2}
\biggr) \Bigl(x_0^k (x_0 + x_2 + \varphi_1 \varphi_2)^k} \Bigr. \\
& & \hspace{1.9in} \Bigl. Y(u,(x_0 + x_2 + \varphi_1 \varphi_2,
\varphi_1))Y(v,(x_2, \varphi_2))w \Bigr)\\ 
&=& x_2^{-1} \delta \biggl(\frac{x_1 - x_0 - \varphi_1 \varphi_2}{x_2}
\biggr) \Bigl(x_0^k (x_0 + x_2 + \varphi_1 \varphi_2)^k \Bigr. \\
& & \hspace{2.4in} \Bigl. Y(Y(u,(x_0, \varphi_1 - \varphi_2))v,(x_2,
\varphi_2))w \Bigr)\\ 
&=& x_2^{-1} \delta \biggl(\frac{x_1 - x_0 - \varphi_1 \varphi_2}{x_2}
\biggr) \Bigl(x_0^k x_1^k Y(Y(u,(x_0, \varphi_1 - \varphi_2))v,(x_2,
\varphi_2))w \Bigr)\\
&=& x_0^k x_1^k x_2^{-1} \delta \biggl(\frac{x_1 - x_0 - \varphi_1
\varphi_2}{x_2} \biggr) Y(Y(u,(x_0, \varphi_1 - \varphi_2))v,(x_2,
\varphi_2))w  \\  
\end{eqnarray*}
which implies the Jacobi identity. $\pfbox$ \\

\section{Expansions of rational superfunctions}

In order to formulate the notions of associativity and
supercommutativity, we will need to interpret correlation functions of
vertex operators with odd formal variables as expansions of certain
rational superfunctions.  In this section we follow and extend the 
treatment of rational functions as presented in \cite{FHL}.

Let $T(U) = \coprod_{n \in \mathbb{N}} T^n(U)$ be the tensor algebra 
over the vector space $U$, where $T^n(U)$ is the $n$-fold tensor 
product of $U$, and let ${\cal J}$ be the ideal of $T(U)$ generated 
by the elements $a \otimes b + b \otimes a$ for $a,b \in U$.  Then 
$\bigwedge (U) = T(U)/{\cal J}$.  (It is understood that $T^0(U) = 
\mathbb{C}$.)  Let $\pi_B$ be the projection {}from $T(U)$ onto $T^0(U)$.  
Then $\pi_B$ is well defined on $\bigwedge(U)$ and is called the 
projection onto the {\it body} of $\bigwedge(U)$ (cf. \cite{D}, 
\cite{B thesis}).  For $a \in \bigwedge_*$, denote $\pi_B (a) = a_B$.

Let $\bigwedge_*[x_1,x_2,...,x_n]_S$ be the ring of rational functions
obtained by inverting (localizing with respect to) the set 
\[S = \biggl\{\sum_{i = 1}^{n} a_i x_i : a_i \in \mbox{$\bigwedge_*^0$}, \; \mbox{not all} \;
(a_i)_B = 0\biggr\} . \] 
Recall the map $\iota_{i_1 ... i_2} :
\mathbb{F}[x_1,...,x_n]_S \longrightarrow \mathbb{F}[[x_1, x_1^{-1},..., x_n,
x_n^{-1}]]$ defined in \cite{FLM} where coefficients of elements in
$S$ are restricted to the field $\mathbb{F}$.  We extend this map to 
$\bigwedge_*[x_1,x_2,...,x_n]_S[\varphi_1,\varphi_2,...,\varphi_n] =
\bigwedge_*[x_1,\varphi_1,x_2,\varphi_2,...,x_n,\varphi_n]_{S}$ in the
obvious way obtaining
\[\iota_{i_1 ... i_2} : \mbox{$\bigwedge_*$} [x_1,\varphi_1,...,x_n,\varphi_n]_S
\longrightarrow \mbox{$\bigwedge_*$}[[x_1, x_1^{-1},..., x_n,
x_n^{-1}]][\varphi_1,...,\varphi_n] .\]
Let $\bigwedge_*[x_1, \varphi_1, x_2, \varphi_2,...,x_n,
\varphi_n]_{S'}$ be the ring of rational functions obtained by
inverting the set   
\[S' = \biggl\{\sum_{\stackrel{i,j = 1}{i<j}}^{n} (a_i x_i + a_{ij} \varphi_i
\varphi_j) : a_i, a_{ij} \in \mbox{$\bigwedge_*^0$}, \; \mbox{not all} \; (a_i)_B
= 0\biggr\}. \]
Since we use the convention that a function of even and odd
variables should be expanded about the even variables, we have 
\[\frac{1}{\sum_{\stackrel{i,j = 1}{i<j}}^{n} (a_i x_i + a_{ij} \varphi_i
\varphi_j)} = \frac{1}{\sum_{i = 1}^{n} a_i x_i} -
\frac{\sum_{\stackrel{i,j = 1}{i<j}}^{n} a_{ij} \varphi_i
\varphi_j}{(\sum_{i = 1}^{n} a_i x_i)^2} .\]
Thus 
\[\mbox{$\bigwedge_*$}[x_1,\varphi_1,x_2,\varphi_2,...,x_n,\varphi_n]_{S'}
\subseteq \mbox{$\bigwedge_*$}[x_1,\varphi_1,x_2,\varphi_2,...,x_n,\varphi_n]_S ,
\]   
and $\iota_{i_1 ... i_2}$ is well defined on
$\bigwedge_*[x_1,\varphi_1,x_2,\varphi_2,...,x_n,\varphi_n]_{S'}$.

In the case $n = 2$, 
\[\iota_{12} : \mbox{$\bigwedge_*$} [x_1,\varphi_1,x_2,\varphi_2]_{S'} 
\longrightarrow \mbox{$\bigwedge_*$} [[x_1,x_2]][\varphi_1,\varphi_2] \]   
is given by first expanding an element of $\bigwedge_*
[x_1,\varphi_1,x_2,\varphi_2]_{S'}$ as a formal series in
$\bigwedge_*[x_1,\varphi_1,x_2,\varphi_2]_S$ and then expanding each
term as a series in $\bigwedge_*[[x_1,x_2]][\varphi_1,\varphi_2]$
containing at most finitely many negative powers of $x_2$ (using
binomial expansions for negative powers of linear polynomials
involving both $x_1$ and $x_2$).

\section{Duality for vertex operator superalgebras}

In \cite{B thesis}, we formulate the notion of $N=1$ supergeometric 
vertex operator superalgebra and show that any such object defines a 
$N=1$ Neveu-Schwarz vertex operator superalgebra with odd formal 
variables.  To show that the alleged vertex operator superalgebra 
satisfies the Jacobi identity, we need the notions of associativity 
and supercommutativity for a vertex operator superalgebra with odd
formal variables.  Together, these notions of associativity and 
(super)commutativity are known as ``duality'', a term which arose 
{}from physics.   Throughout this section we follow and extend the 
treatment of duality as presented in \cite{FHL}.

Let $(V, Y(\cdot,(x,\varphi)), \mathbf{1}, \tau)$ be a vertex operator
algebra with odd formal variables.  Let $V_{(n)}^*$ be the dual module
of $V_{(n)}$ for $n \in \frac{1}{2} \mathbb{Z}$, i.e.,$V_{(n)}^* =
\mbox{Hom}_{\bigwedge_*} (V, \bigwedge_*)$.  Let  
\[V' = \coprod_{n \in \frac{1}{2} \mathbb{Z}} V_{(n)}^* \]
be the graded dual space of $V$, 
\[\bar{V} = \prod _{n \in \frac{1}{2} \mathbb{Z}} V_{(n)} = V'^* \] 
the algebraic completion of $V$, and $\langle \cdot , \cdot \rangle$
the natural pairing between $V'$ and $\bar{V}$.  We now formulate the
weak supercommutativity and weak associativity properties of a vertex
operator superalgebra with odd formal variables into slightly stronger
statements about ``matrix coefficients'' of products and iterates of
vertex operators with odd formal variables.  

\begin{prop}\label{supercommutativity} 
{\bf (a) (rationality of products)}  For $u, v, w \in V$, with $u$, and $v$
of homogeneous sign, and $v' \in V'$, the formal series 
\[ \langle v', Y(u,(x_1,\varphi_1)) Y(v,(x_2,\varphi_2)) w \rangle, \]
which involves
only finitely many negative powers of $x_2$ and only finitely many
positive powers of $x_1$, lies in the image of the map $\iota_{12}$: 
\[\langle v', Y(u,(x_1,\varphi_1)) Y(v,(x_2,\varphi_2)) w \rangle
= \iota_{12} f(x_1,\varphi_1,x_2,\varphi_2) , \]
where the (uniquely determined) element $f \in \bigwedge_*
[x_1,\varphi_1,x_2,\varphi_2]_{S'}$ is of the form  
\[f(x_1,\varphi_1,x_2,\varphi_2) =
\frac{g(x_1,\varphi_1,x_2,\varphi_2)}{x_1^r x_2^s (x_1 - x_2 - 
\varphi_1 \varphi_2)^t} \]
for some $g \in \bigwedge_* [x_1,\varphi_1,x_2,\varphi_2]$ and $r, s, t
\in \mathbb{Z}$.  

{\bf (b) (supercommutativity)}  We also have
\[\langle v', Y(v,(x_2,\varphi_2))Y(u,(x_1,\varphi_1)) w \rangle =
(-1)^{\eta(u)\eta(v)} \iota_{21} f(x_1,\varphi_1,x_2,\varphi_2)
, \] 
i.e,
\[\iota_{12}^{-1} \langle v', Y(u,(x_1,\varphi_1))
Y(v,(x_2,\varphi_2)) w \rangle = \hspace{2.3in} \]
\[\hspace{1.6in} (-1)^{\eta(u)\eta(v)} \iota_{21}^{-1} \langle v',
Y(v,(x_2,\varphi_2)) Y(u,(x_1,\varphi_1)) w \rangle .\]  
\end{prop}

{\it Proof:} \hspace{.2cm}  Part (a) follows {}from the positive
energy axiom (\ref{positive energy}) and truncation condition
(\ref{truncation}) for a vertex operator superalgebra.   
For part (b), we note that by weak supercommutativity, there 
exists $k \in \Z$ such that
\begin{equation}\label{commute}
(x_1 - x_2 - \varphi_1 \varphi_2)^k \langle v', Y(u,(x_1,\varphi_1))
Y(v,(x_2,\varphi_2)) w\rangle = \hspace{1in}  
\end{equation}
\[ \hspace{.7in} (-1)^{\eta(u)\eta(v)} (x_1 - x_2 - \varphi_1
\varphi_2)^k \langle v', Y(v,(x_2,\varphi_2)) Y(u,(x_1,\varphi_1)) w
\rangle \] 
for all $w \in V$ and $v' \in V'$. {}From (a), we know the left-hand side of (\ref{commute}) involves only
finitely many negative powers of $x_2$ and only finitely many positive
powers of $x_1$.  However, the right-hand side of (\ref{commute})
involves only finitely many negative powers of $x_1$ and only finitely
many positive powers of $x_2$.  Thus multiplying both sides of
(\ref{commute}) by $(x_1 - x_2 - \varphi_1 \varphi_2)^{-k}$ results in
well-defined power series as long as on the left-hand side we expand 
$(x_1 - x_2 - \varphi_1 \varphi_2)^{-k}$ in positive powers of $x_2$ 
and on the right-hand side we expand $(x_1 - x_2 - \varphi_1 \varphi_2)^{-k}$ in
positive powers of $x_1$.  The result follows. $\pfbox$ \\

\begin{prop} \label{iterates}
{\bf (a) (rationality of iterates)} For $u, v, w \in V$, and $v' \in
V'$, the formal series $\langle v', Y(Y(u,(x_0,\varphi_1 -
\varphi_2))v,(x_2,\varphi_2)) w \rangle$, which involves only finitely
many negative powers of $x_0$ and only finitely many positive powers
of $x_2$, lies in the image of the map $\iota_{20}$: 
\[\langle v', Y(Y(u,(x_0,\varphi_1 - \varphi_2))v,(x_2,\varphi_2))w
\rangle  = \iota_{20} h(x_0,\varphi_1 - \varphi_2,x_2,\varphi_2) , \]
where the (uniquely determined) element $h \in \bigwedge_*
[x_0,\varphi_1,x_2,\varphi_2]_{S'}$ is of the form      
\[h(x_0,\varphi_1 - \varphi_2,x_2,\varphi_2) = \frac{k(x_0,\varphi_1 -
\varphi_2,x_2,\varphi_2)}{x_0^r x_2^s (x_0 + x_2 - \varphi_1
\varphi_2)^t} \]   
for some $k \in \bigwedge_* [x_0,\varphi_1,x_2,\varphi_2]$ and $r, s, t
\in \mathbb{Z}$.  

{\bf (b)}  The formal series $\langle v', Y(u,(x_0 + x_2  + \varphi_1
\varphi_2, \varphi_1)) Y(v,(x_2,\varphi_2)) w \rangle$, 
which involves only finitely many negative powers of $x_2$ and only
finitely many positive powers of $x_0$, lies in the image of
$\iota_{02}$, and in fact  
\[\langle v', Y(u,(x_0 + x_2  + \varphi_1 \varphi_2, \varphi_1))
Y(v,(x_2,\varphi_2)) w \rangle =  \iota_{02} h(x_0,\varphi_1 -
\varphi_2, x_2, \varphi_2) . \]    
\end{prop}

{\it Proof:} \hspace{.2cm}  Part (a) follows {}from the positive energy
axiom (\ref{positive energy}) and truncation condition
(\ref{truncation}) for a vertex operator superalgebra. 
For part (b), we note that {}from weak associativity, there exists $k \in \Z$ such that
\begin{equation}\label{assoc}
(x_0 + x_2 + \varphi_1 \varphi_2)^k \langle v', Y(Y(u,(x_0,\varphi_1 -
\varphi_2))v,(x_2,\varphi_2))w \rangle = \hspace{1in}  
\end{equation}
\[ \hspace{.7in} (x_0 + x_2 + \varphi_1 \varphi_2)^k \langle v',
Y(u,(x_0 + x_2  + \varphi_1 \varphi_2, \varphi_1))
Y(v,(x_2,\varphi_2)) w  \rangle \] 
for all $v' \in V'$. {}From (a), we know the left-hand side of
(\ref{assoc}) involves only finitely many negative powers of $x_0$ and
only finitely many positive powers of $x_2$.  However, the right-hand
side of (\ref{assoc}) involves only finitely many negative powers of
$x_2$ and only finitely many positive powers of $x_0$.  Thus
multiplying both sides of (\ref{commute}) by $(x_0 + x_2 + \varphi_1
\varphi_2)^{-k}$ results in a well-defined power series as long as on the left-hand side we
expand $(x_0 + x_2 + \varphi_1 \varphi_2)^{-k}$ in positive powers of
$x_0$ and on the right-hand side we expand $(x_0 + x_2 + \varphi_1
\varphi_2)^{-k}$ in positive powers of $x_2$.  The result follows.
$\pfbox$ \\    

\begin{prop}\label{associativity} 
{\bf (associativity)}  We have the following equality of rational
functions: 
\[\iota^{-1}_{12} \langle v', Y(u,(x_1,\varphi_1))
Y(v,(x_2,\varphi_2)) w \rangle = \hspace{2.4in} \]
\[\hspace{1in} \left. \left( \iota^{-1}_{20} \langle v',
Y(Y(u,(x_0,\varphi_1 - \varphi_2))v,(x_2,\varphi_2)) w \rangle \right)
\right|_{x_0 = x_1 - x_2 - \varphi_1 \varphi_2} \] 
\end{prop}

{\it Proof:} \hspace{.2cm}  Let $f(x_1,\varphi_1,x_2,\varphi_2)$ be
the rational function in Proposition \ref{supercommutativity}.  Then
$f$ satisfies 
\[\iota_{02} f(x_0 + x_2 + \varphi_1
\varphi_2,\varphi_1,x_2,\varphi_2) = \left. \left( \iota_{12}
f(x_1,\varphi_1,x_2,\varphi_2) \right) \right|_{x_1 = x_0 + x_2 +
\varphi_1 \varphi_2} . \] 
Thus for $h(x_0,\varphi_1 -
\varphi_2, x_2, \varphi_2)$ {}from Proposition \ref{iterates}, we have
$h(x_0,\varphi_1 - \varphi_2, x_2, \varphi_2) = f(x_0 + x_2 +
\varphi_1 \varphi_2,\varphi_1,x_2,\varphi_2)$.  The result
follows {}from Propositions \ref{supercommutativity} and
\ref{iterates}. $\pfbox$ \\  

Note that rationality of products and iterates and supercommutativity
and associativity imply weak supercommutativity and weak
associativity.   Thus by Proposition \ref{equals Jacobi1}, we have: 

\begin{prop}\label{duality}
In the presence of the other axioms in the definition of
vertex operator superalgebra with odd variables, the Jacobi identity
follows {}from the rationality of products and iterates and
supercommutativity and associativity.  In particular, the Jacobi
identity may be replaced by these properties. 
\end{prop}

This can also be proved by using the delta-function identity
(\ref{delta 3 terms with phis}) and the substitution rule (\ref{delta
substitute}).

\end{document}